   \documentstyle[12pt]{article}
    \title{{\bf  A theory of tensor products for module categories for
a vertex operator algebra, III}}
    \author{Yi-Zhi Huang and James Lepowsky}
    \date{}
    
\begin{document}
    \bibliographystyle{alpha}
    \maketitle

    \input amssym.def
    \input amssym
    \newtheorem{rema}{Remark}[section]
    \newtheorem{propo}[rema]{Proposition}
    \newtheorem{theo}[rema]{Theorem}
   \newtheorem{defi}[rema]{Definition}
    \newtheorem{lemma}[rema]{Lemma}
    \newtheorem{corol}[rema]{Corollary}
     \newtheorem{exam}[rema]{Example}
	\newcommand{\ba}{\begin{array}}
	\newcommand{\ea}{\end{array}}
        \newcommand{\be}{\begin{equation}}
        \newcommand{\ee}{\end{equation}}
	\newcommand{\bea}{\begin{eqnarray}}
	\newcommand{\eea}{\end{eqnarray}}
	\newcommand{\nno}{\nonumber}
	\newcommand{\lbar}{\bigg\vert}
	\newcommand{\p}{\partial}
	\newcommand{\dps}{\displaystyle}
	\newcommand{\bra}{\langle}
	\newcommand{\ket}{\rangle}
 \newcommand{\res}{\mbox{ \rm Res}}
\newcommand{\wt}{\mbox{ \rm wt}\ }
 \newcommand{\pf}{{\it Proof}\hspace{2ex}}
 \newcommand{\epf}{\hspace{2em}$\Box$}
 \newcommand{\epfv}{\hspace{1em}$\Box$\vspace{1em}}

\hyphenation{Phil-a-del-phia}

 \makeatletter
\newlength{\@pxlwd} \newlength{\@rulewd} \newlength{\@pxlht}
\catcode`.=\active \catcode`B=\active \catcode`:=\active \catcode`|=\active
\def\sprite#1(#2,#3)[#4,#5]{
   \edef\@sprbox{\expandafter\@cdr\string#1\@nil @box}
   \expandafter\newsavebox\csname\@sprbox\endcsname
   \edef#1{\expandafter\usebox\csname\@sprbox\endcsname}
   \expandafter\setbox\csname\@sprbox\endcsname =\hbox\bgroup
   \vbox\bgroup
  \catcode`.=\active\catcode`B=\active\catcode`:=\active\catcode`|=\active
      \@pxlwd=#4 \divide\@pxlwd by #3 \@rulewd=\@pxlwd
      \@pxlht=#5 \divide\@pxlht by #2
      \def .{\hskip \@pxlwd \ignorespaces}
      \def B{\@ifnextchar B{\advance\@rulewd by \@pxlwd}{\vrule
         height \@pxlht width \@rulewd depth 0 pt \@rulewd=\@pxlwd}}
      \def :{\hbox\bgroup\vrule height \@pxlht width 0pt depth
0pt\ignorespaces}
      \def |{\vrule height \@pxlht width 0pt depth 0pt\egroup
         \prevdepth= -1000 pt}
   }
\def\endsprite{\egroup\egroup}
\catcode`.=12 \catcode`B=11 \catcode`:=12 \catcode`|=12\relax
\makeatother

\def\hboxtr{\FormOfHboxtr} 
\sprite{\FormOfHboxtr}(25,25)[0.5 em, 1.2 ex] 

:BBBBBBBBBBBBBBBBBBBBBBBBB |
:BB......................B |
:B.B.....................B |
:B..B....................B |
:B...B...................B |
:B....B..................B |
:B.....B.................B |
:B......B................B |
:B.......B...............B |
:B........B..............B |
:B.........B.............B |
:B..........B............B |
:B...........B...........B |
:B............B..........B |
:B.............B.........B |
:B..............B........B |
:B...............B.......B |
:B................B......B |
:B.................B.....B |
:B..................B....B |
:B...................B...B |
:B....................B..B |
:B.....................B.B |
:B......................BB |
:BBBBBBBBBBBBBBBBBBBBBBBBB |

\endsprite

\def\shboxtr{\FormOfShboxtr} 
\sprite{\FormOfShboxtr}(25,25)[0.3 em, 0.72 ex] 

:BBBBBBBBBBBBBBBBBBBBBBBBB |
:BB......................B |
:B.B.....................B |
:B..B....................B |
:B...B...................B |
:B....B..................B |
:B.....B.................B |
:B......B................B |
:B.......B...............B |
:B........B..............B |
:B.........B.............B |
:B..........B............B |
:B...........B...........B |
:B............B..........B |
:B.............B.........B |
:B..............B........B |
:B...............B.......B |
:B................B......B |
:B.................B.....B |
:B..................B....B |
:B...................B...B |
:B....................B..B |
:B.....................B.B |
:B......................BB |
:BBBBBBBBBBBBBBBBBBBBBBBBB |

\endsprite

\vspace{2em}

\begin{abstract}
This is the third part in a series of papers developing a tensor
product theory for modules for a vertex operator algebra.  The goal of
this theory is to construct a ``vertex tensor category'' structure on
the category of modules for a suitable vertex operator algebra.  The
notion of vertex tensor category is essentially a ``complex analogue''
of the notion of symmetric tensor category, and in fact a vertex
tensor category produces a braided tensor category in a natural way.
In this paper, we focus on a particular element $P(z)$ of a certain
moduli space of three-punctured Riemann spheres; in general, every
element of this moduli space will give rise to a notion of tensor
product, and one must consider all these notions in order to construct
a vertex tensor category.  Here we present the fundamental properties
of the $P(z)$-tensor product of two modules for a vertex operator
algebra.  We give two constructions of a $P(z)$-tensor product, using
the results, established in Parts I and II of this series,
for a certain other element of the moduli space.  The definitions and
results in Parts I and II are recalled.
\end{abstract}

\vspace{2em}

\renewcommand{\theequation}{I.\arabic{equation}}
\renewcommand{\therema}{I.\arabic{rema}}
\setcounter{equation}{0}
\setcounter{rema}{0}

One of the most important operations in the representation theory of
Lie algebras is the tensor product operation for modules for a Lie
algebra. Together with this operation and the coefficient
field as an identity object, the category of modules for a Lie algebra
is a symmetric tensor category. But for any category of modules
 of a fixed nonzero level for an
affine Lie algebra, the usual tensor product
operation  does not give a tensor category
structure. Instead, there are certain categories of modules of a fixed
level for an affine Lie algebra equipped with
conformal-field-theoretic tensor product operations, and modules
playing the role of identity objects, giving braided tensor
categories.  {}From the viewpoint of conformal field theory,
the most relevant cases involve
positive integral levels, including
the case of the
category whose objects are finite direct sums of modules isomorphic to
standard (integrable highest weight) modules of a
fixed positive integral level for an affine Lie algebra.
This was explained on a physical level of rigor by Moore and Seiberg
\cite{MS}, under the very subtle and nontrivial
assumption that the ``chiral vertex operators'' have suitable
``operator product expansion'' properties, or that an equivalent geometric
axiom in conformal field theory holds. In \cite{KL}--\cite{KL4},
Kazhdan and Lusztig constructed the braided tensor category structure
for another category of modules, of a fixed but sufficiently negative
level, for an affine Lie algebra.

Vertex operator algebras (\cite{B}, \cite{FLM}, \cite{FHL}) are
analogous to both Lie algebras and commutative associative algebras.
They are essentially equivalent to ``chiral algebras'' in conformal
field theory (see for example \cite{MS}).  Motivated partly by the
analogy between vertex operator algebras and Lie algebras and partly
by the announcement \cite{KL}, we initiated a theory of tensor
products for modules for a vertex operator algebra in
\cite{HL0}.  This theory was developed in detail
beginning in \cite{HL1}--\cite{HL2}, and an overview of the theory was
given in \cite{HL3}, where it was announced that for a vertex operator
algebra satisfying suitable conditions, its module category has a
natural structure of ``vertex tensor category.''  It was also
announced there that the underlying category of a vertex tensor
category has a natural structure of braided tensor category.  In
particular, the category of modules for the vertex operator algebra
associated to a minimal model and the category whose objects are
finite direct sums of modules isomorphic to standard modules of a
fixed positive integral level for an affine Lie algebra have natural
braided tensor category structure. In this theory, instead of being an
assumption as in \cite{MS}, the operator product expansion (or
associativity) for chiral vertex operators (or intertwining operators)
is a consequence.

The present paper (Part III) is a continuation of \cite{HL1} (Part I) and
\cite{HL2} (Part II), to which---especially \cite{HL1}---the
reader is referred for the necessary background, including references.
The reader is referred to \cite{HL3}
for the motivation and  description of the main
results of the tensor product theory developed by the present series of papers.
To make the present paper as self-contained as reasonably possible, we
shall provide in this introduction a systematic summary of the necessary
results of Parts I and II.

In Part I, the
notions of $P(z)$- and $Q(z)$-tensor product ($z\in {\Bbb C}^{\times}$)
of modules for a vertex
operator algebra were introduced, and under suitable conditions, two
constructions
of a $Q(z)$-tensor product were given in Part I based on certain results
proved in Part II.
(The symbols $P(z)$ and $Q(z)$ designate two elements of a certain
moduli space of spheres with punctures and local coordinates.)
In the present paper, the notion of $P(z)$-tensor product
is discussed in parallel with the discussion in Section 4 of Part I for
that of $Q(z)$-tensor product,   and two constructions
of a $P(z)$-tensor product are given,
using the results for the $Q(z^{-1})$-tensor product.

We now describe some of our basic notation and elementary tools.
We work over ${\Bbb C}$.
In this paper, as in \cite{HL1} and \cite{HL2},
the symbols $x, x_{0}, x_{1}, \dots$ and $t$ are independent commuting
formal variables, and all expressions involving these variables are
to be understood as formal Laurent series. We use the ``formal
$\delta$-function''
$$\delta(x)=\sum_{n\in {\Bbb Z}}x^{n}.$$
It has the following simple and
fundamental property:
For any $f(x)\in {\Bbb C}[x, x^{-1}]$,
$$f(x)\delta(x)=f(1)\delta(x).$$
This property has many important variants. For example, for any
$$X(x_{1},
x_{2})\in (\mbox{End }W)[[x_{1}, x_{1}^{-1}, x_{2}, x_{2}^{-1}]]$$
(where $W$ is a vector space) such that
$$
\lim_{x_{1}\to x_{2}}X(x_{1}, x_{2})=X(x_{1},
x_{2})\lbar_{x_{1}=x_{2}}
$$
exists, we have
$$
X(x_{1}, x_{2})\delta\left(\frac{x_{1}}{x_{2}}\right)=X(x_{2}, x_{2})
\delta\left(\frac{x_{1}}{x_{2}}\right).
$$
The existence of this ``algebraic limit'' means that
for an arbitrary vector $w\in W$, the coefficient of each power of
$x_{2}$ in the formal expansion $X(x_{1}, x_{2})w\lbar_{x_{1}=x_{2}}$
is a finite sum.  We use the convention that negative powers of a
binomial are to be expanded in nonnegative powers of the second
summand. For example,
$$
x_{0}^{-1}\delta\left(\frac{x_{1}-x_{2}}{x_{0}}\right)=\sum_{n\in {\Bbb Z}}
\frac{(x_{1}-x_{2})^{n}}{x_{0}^{n+1}}=\sum_{m\in {\Bbb N},\; n\in {\Bbb Z}}
(-1)^{m}{{n}\choose {m}} x_{0}^{-n-1}x_{1}^{n-m}x_{2}^{m}.
$$
We have the following identities:
\begin{eqnarray*}
&{\dps x_{1}^{-1}\delta\left(\frac{x_{2}+x_{0}}{x_{1}}\right)=x_{2}^{-1}\left(
\frac{x_{1}-x_{0}}{x_{2}}\right),}&\\
&{\dps x_{0}^{-1}\delta\left(\frac{x_{1}-x_{2}}{x_{0}}\right)-
x_{0}^{-1}\delta\left(\frac{x_{2}-x_{1}}{-x_{0}}\right)=
x_{2}^{-1}\delta\left(\frac{x_{1}-x_{0}}{x_{2}}\right).}&
\end{eqnarray*}
We shall use these properties and identities later on
without explicit comment.  Here and below, it is important to note
that the
relevant sums and products, etc., of formal series, are well defined.
See \cite{FLM}, \cite{FHL}, \cite{HL1} and
\cite{HL2} for further
discussion and many examples of their use, including their role in formulating
and using the Jacobi identity for vertex operator algebras, modules and
intertwining operators.

As in \cite{HL1} and \cite{HL2}, the symbol $z$ will always
denote a nonzero complex number. We shall always choose $\log z$ such that
$$\log z=\log |z| +i\arg z \;\;\;\;\mbox{with}\;\;\;\;0\le \arg z<2\pi.$$
Arbitrary values of the log function at $z$ will be denoted
$$l_{p}(z)=\log z+2p\pi i$$
for $p\in {\Bbb Z}$.

We fix a vertex operator algebra $V \; (= (V, Y, {\bf 1}, \omega))$;
recall that $Y$ is the vertex operator map; for $v\in V$, $Y(v,
x)=\sum_{n\in {\Bbb Z}}v_{n}x^{-n-1}$ ($v_{n}\in \mbox{\rm End}\ V$);
${\bf 1}$ is the vacuum vector; and $\omega$ is the element whose
vertex operator $Y(\omega,x)=\sum_{n\in {\Bbb Z}}L(n)x^{-n-2}$ gives
the Virasoro algebra.  (See
\cite{FLM} and \cite{FHL} for the basic definitions.)  The symbols $W,
W_{1}, W_{2}, \dots$ will denote (${\Bbb C}$-graded) $V$-modules.  The
symbol $Y$ will denote the vertex operator map for $W$ (as well as for
$V$), and the symbols $Y_{1}, Y_{2}, \dots$ will denote the vertex
operator maps for the $V$-modules $W_{1}, W_{2}, \dots$, respectively.
When necessary, we shall use notation such as $(W, Y)$ to designate a
$V$-module.  The symbols ${\cal Y}, {\cal Y}_{1},
\dots$ will denote intertwining operators (as defined in \cite{FHL}).
The vector space of intertwining operators of type ${W_{3}}\choose
{W_{1}W_{2}}$ will be denoted ${\cal V}^{W_{3}}_{W_{1}W_{2}}$.  The
dimension of ${\cal V}^{W_{3}}_{W_{1}W_{2}}$ is the corresponding {\it
fusion rule}.  As in \cite{HL1}, for any $v\in V$,
$$Y_{t}(v, x)=v\otimes t^{-1}\delta\left(\frac{x}{t}\right)\in
V\otimes {\Bbb C}[[t, t^{-1}, x, x^{-1}]].$$
We also need the notion of generalized
module for $V$. A generalized module for $V$ is a pair $(W, Y)$
satisfying all the axioms for a $V$-module except for the two axioms
on the grading of $W$ (see Definition 2.11 in \cite{HL1}).

Now we present the main results on $Q(z)$-tensor products, together
with the necessary definitions and explanations, {}from Parts I and II.

Fix a nonzero complex number $z$ and let $(W_{1}, Y_{1})$ and $(W_{2},
Y_{2})$ be $V$-modules. We recall the definition of $Q(z)$-tensor product
of $W_{1}$ and $W_{2}$.
A {\it $Q(z)$-intertwining map of type ${W_{3}}\choose {W_{1}W_{2}}$}
is a linear map
$F: W_{1}\otimes W_{2} \to \overline{W}_{3}$ ($\overline{W}$ being the
formal algebraic completion of a given module $W$) such that
\begin{eqnarray*}
\lefteqn{z^{-1}\delta\left(\frac{x_{1}-x_{0}}{z}\right)
Y^{*}_{3}(v, x_{0})F(w_{(1)}\otimes w_{(2)})=}\nonumber\\
&&=x_{0}^{-1}\delta\left(\frac{x_{1}-z}{x_{0}}\right)
F(Y_{1}^{*}(v, x_{1})w_{(1)}\otimes w_{(2)})\nonumber\\
&&\hspace{2em}-x_{0}^{-1}\delta\left(\frac{z-x_{1}}{-x_{0}}\right)
F(w_{(1)}\otimes Y_{2}(v, x_{1})w_{(2)})
\end{eqnarray*}
for $v\in V$, $w_{(1)}\in W_{1}$, $w_{(2)}\in W_{2}$, where for a
vertex operator $Y$,
$$Y^*(v,x) = Y(e^{xL(1)}(-x^{-2})^{L(0)}v, x^{-1}).$$
We denote the vector space of $Q(z)$-intertwining maps of type
${W_{3}}\choose {W_{1}W_{2}}$ by ${\cal M}[Q(z)]^{W_{3}}_{W_{1}W_{2}}$.

We define a {\it $Q(z)$-product of $W_{1}$ and $W_{2}$} to be a
$V$-module $(W_{3}, Y_{3})$ together with a $Q(z)$-intertwining map
$F$ of type ${W_{3}}\choose {W_{1}W_{2}}$ and we denote it by $(W_{3},
Y_{3}; F)$ (or $(W_{3}, F)$).
 Let $(W_{3}, Y_{3}; F)$ and $(W_{4}, Y_{4}; G)$ be two
$Q(z)$-products of $W_{1}$ and $W_{2}$. A {\it morphism} {from}
$(W_{3}, Y_{3}; F)$ to $(W_{4}, Y_{4}; G)$ is a module map $\eta$ {from}
$W_{3}$ to $W_{4}$ such that
$$G=\overline{\eta}\circ F$$
where $\overline{\eta}$ is the natural map {from} $\overline{W}_{3}$ to
$\overline{W}_{4}$ uniquely extending $\eta$.
A {\it $Q(z)$-tensor product of $W_{1}$ and $W_{2}$} is
a $Q(z)$-product $(W_{1}\boxtimes_{Q(z)} W_{2},
Y_{Q(z)}; \boxtimes_{Q(z)})$
such that for any $Q(z)$-product
$(W_{3}, Y_{3}; F)$, there is a unique morphism {from}
$(W_{1}\boxtimes_{Q(z)} W_{2}, Y_{Q(z)};
\boxtimes_{Q(z)})$ to $(W_{3}, Y_{3}; F)$.
The $V$-module $(W_{1}\boxtimes_{Q(z)} W_{2},  Y_{Q(z)})$ is
called a {\it $Q(z)$-tensor product module} of $W_{1}$ and $W_{2}$.
If it exists, it is unique up to canonical isomorphism.

We now describe the close connection between intertwining operators of
type ${W'_{1}}\choose {W'_{3}W_{2}}$ and $Q(z)$-intertwining maps of
type ${W_{3}}\choose {W_{1}W_{2}}$, where $W'$ denotes the contragredient
of a module $W$, as defined in \cite{FHL} (or, more generally, the
graded dual space of a graded vector space).
Fix an integer $p$. Let ${\cal
Y}$ be an intertwining operator of type ${W'_{1}}\choose
{W'_{3}W_{2}}$.  We have a linear map $F^{Q(z)}_{{\cal Y}, p}:
W_{1}\otimes W_{2}\to \overline{W}_{3}$ determined by the condition
\begin{equation}\label{4.9}
\langle w'_{(3)},
F^{Q(z)}_{{\cal Y}, p}(w_{(1)}\otimes w_{(2)})\rangle_{W_{3}} =\langle
w_{(1)}, {\cal Y}(w'_{(3)}, e^{l_{p}(z)})w_{(2)}\rangle_{W'_{1}}
\end{equation}
for all $w_{(1)}\in W_{1}$, $w_{(2)}\in W_{2}$, $w'_{(3)}\in W'_{3}$, where
for a module $W$, $\langle\cdot, \cdot\rangle_{W}$ denotes the
canonical pairing
between $W'$ and $\overline{W}$, and where the notation
${\cal Y}(\cdot, e^{\zeta})$ for $\zeta\in {\Bbb C}$ is shorthand for
${\cal Y}(\cdot, x)|_{x^{n}=e^{n\zeta},\;n\in {\Bbb C}}$, which is
well defined; note that ${\cal Y}(\cdot, e^{\zeta})$ actually depends
on $\zeta$ and not just $e^{\zeta}$.
Using the Jacobi identity for
${\cal Y}$, we see easily that $F^{Q(z)}_{{\cal Y}, p}$ is a
$Q(z)$-intertwining map of
type ${W_{3}}\choose {W_{1}W_{2}}$.

Conversely, given
 a $Q(z)$-intertwining map $F$ of type ${W_{3}}\choose {W_{1}W_{2}}$,
as a linear map {from} $W_{1}\otimes W_{2}$ to $\overline{W}_{3}$ it
gives us an element of $(W_{1}\otimes W'_{3}\otimes W_{2})^{*}$ whose
value at $w_{(1)}\otimes w'_{(3)}\otimes w_{(2)}$ is $$\langle
w'_{(3)}, F(w_{(1)}\otimes w_{(2)})\rangle_{W_3}.$$ But since every
element of $(W_{1}\otimes W'_{3}\otimes W_{2})^{*}$ also amounts to a
linear map {from} $W'_{3}\otimes W_{2}$ to $W^{*}_{1}$, we have such a
map as well. Let $w'_{(3)}\in W'_{3}$ and $w_{(2)}\in W_{2}$ be
homogeneous elements. Since $W^{*}_{1}=\prod_{n\in {\Bbb
C}}(W'_{1})_{(n)}$, the image of $w'_{(3)}\otimes w_{(2)}$ under our
map can be written as $\sum_{n\in {\Bbb
C}}(w'_{(3)})_{n}w_{(2)}e^{(-n-1)l_{p}(z)}$ where for any $n\in {\Bbb
C}$, $(w'_{(3)})_{n}w_{(2)}e^{(-n-1)l_{p}(z)}$ is the projection of
the image to the homogeneous subspace of $W'_{1}$ of weight equal to
$$\mbox{\rm wt}\;w'_{(3)}-n-1+\mbox{\rm wt}\; w_{(2)}.$$ (Here we are
defining elements denoted $(w'_{(3)})_{n}w_{(2)}$ of $W'_{1}$ for
$n\in {\Bbb C}$.)  We define
$${\cal Y}_{F,p}(w'_{(3)}, x)w_{(2)}
=\sum_{n\in {\Bbb C}}(w'_{(3)})_{n}w_{(2)}x^{-n-1}\in W'_{1}\{ x\}$$
for all homogeneous elements $w'_{(3)}\in W'_{3}$ and $w_{(2)}\in
W_{2}$.
(For a vector space $W$, we use the notation
$$W\{ x\} = \{\sum_{n\in {\Bbb C}}a_n x^n \; | \; a_n \in W, \;
n \in {\Bbb C} \};$$
in particular, we are allowing complex powers of our commuting formal
variables.)
Using linearity, we extend ${\cal Y}_{F,p}$ to a linear map
\begin{eqnarray*}
W'_{3}\otimes W_{2}&\to& W'_{1}\{ x\}\nno\\
w'_{(3)}\otimes w_{(2)}&\mapsto &{\cal Y}_{F,p}(w'_{(3)}, x)w_{(2)}.
\end{eqnarray*}
The correspondence $F\mapsto {\cal Y}_{F, p}$ is linear, and
we have ${\cal Y}_{F^{Q(z)}_{{\cal Y}, p}, p}={\cal Y}$
for an intertwining operator
${\cal Y}$ of
type ${W'_{1}}\choose {W'_{3}W_{2}}$.

\begin{propo}\label{4-7}
For $p\in {\Bbb Z}$, the correspondence ${\cal Y}\mapsto
F^{Q(z)}_{{\cal Y}, p}$ is
a linear isomorphism {from} the space ${\cal
V}^{W'_{1}}_{W'_{3}W_{2}}$ of intertwining operators of type
${W'_{1}}\choose {W'_{3}\; W_{2}}$ to the space
${\cal M}[Q(z)]^{W_{3}}_{W_{1}W_{2}}$ of
$Q(z)$-intertwining maps of type ${W_{3}}\choose {W_{1}W_{2}}$. Its
inverse is given by $F\mapsto {\cal Y}_{F, p}$.
\end{propo}

The following immediate result relates module maps {from} a tensor
product module with intertwining maps and intertwining operators:
\begin{propo}
Suppose that $W_{1}\boxtimes_{Q(z)}W_{2}$ exists. We have a natural
isomorphism
\begin{eqnarray*}
\mbox{\rm Hom}_{V}(W_{1}\boxtimes_{Q(z)}W_{2}, W_{3})&\stackrel{\sim}{\to}&
{\cal M}^{W_{3}}_{W_{1}W_{2}}\nno\\
\eta&\mapsto& \overline{\eta}\circ \boxtimes_{Q(z)}
\end{eqnarray*}
and for $p\in {\Bbb Z}$, a natural isomorphism
\begin{eqnarray*}
\mbox{\rm Hom}_{V}(W_{1}\boxtimes_{Q(z)} W_{2}, W_{3})
&\stackrel{\sim}{\rightarrow}&
{\cal V}^{W'_{1}}_{W'_{3}W_{2}}\nno\\
\eta&\mapsto & {\cal Y}_{\eta, p}
\end{eqnarray*}
where ${\cal Y}_{\eta, p}={\cal Y}_{F, p}$ with
$F=\overline{\eta}\circ \boxtimes_{Q(z)}$.
\end{propo}

We have:

\begin{propo}\label{4-9}
For any integer $r$, there is a natural isomorphism
$$B_{r}: {\cal V}^{W_{3}}_{W_{1}W_{2}}\to
{\cal V}^{W'_{1}}_{W'_{3}W_{2}}$$
defined by the condition  that for any
intertwining operator ${\cal Y}$ in ${\cal V}^{W_{3}}_{W_{1}W_{2}}$ and
$w_{(1)}\in W_{1}$, $w_{(2)}\in W_{2}$, $w'_{(3)}\in W'_{3}$,
\begin{eqnarray*}
\lefteqn{\langle w_{(1)}, B_{r}({\cal Y})(w'_{(3)}, x)
w_{(2)}\rangle_{W'_1}=}\nno\\
&&=\langle e^{-x^{-1}L(1)}w'_{(3)}, {\cal Y}(e^{xL(1)}w_{(1)}, x^{-1})
\cdot\nno\\
&&\hspace{4em}\cdot e^{-xL(1)}e^{(2r+1)\pi iL(0)}
x^{-2L(0)}w_{(2)}\rangle_{W_3}.
\end{eqnarray*}
\end{propo}

The last two results give:

\begin{corol}
 For any $V$-modules $W_{1}$, $W_{2}$, $W_{3}$ such that
$W_{1}\boxtimes_{Q(z)}W_{2}$ exists and any integers $p$ and $r$,
we have a natural isomorphism
\begin{eqnarray*}
\mbox{\rm Hom}_{V}(W_{1}\boxtimes_{Q(z)} W_{2}, W_{3})&
\stackrel{\sim}{\rightarrow}&
{\cal V}^{W_{3}}_{W_{1}W_{2}}\nno\\
\eta&\mapsto &B^{-1}_{r}({\cal Y}_{\eta, p}).
\end{eqnarray*}
\end{corol}

It is clear {from} the definition of $Q(z)$-tensor product that
the $Q(z)$-tensor product operation
distributes over direct sums in the following sense:

\begin{propo}
For $V$-modules $U_{1}, \dots, U_{k}$, $W_{1}, \dots, W_{l}$, suppose
that each $U_{i}\boxtimes_{Q(z)}W_{j}$ exists. Then
$(\coprod_{i}U_{i})\boxtimes_{Q(z)}(\coprod_{j}W_{j})$ exists and there is a
natural isomorphism
$$\biggl(\coprod_{i}U_{i}\biggr)\boxtimes_{Q(z)}
\biggl(\coprod_{j}W_{j}\biggr)\stackrel{\sim}
{\rightarrow} \coprod_{i,j}U_{i}\boxtimes_{Q(z)}W_{j}.\hspace{2em} $$
\end{propo}

Now consider $V$-modules $W_{1}$, $W_{2}$ and $W_{3}$ and suppose that
$$\dim {\cal M}[Q(z)]_{W_{1}W_{2}}^{W_{3}}<\infty.$$
The natural
evaluation map
\begin{eqnarray*}
W_{1}\otimes W_{2}\otimes {\cal M}[Q(z)]^{W_{3}}_{W_{1}W_{2}}&\to&
\overline{W}_{3}
\nno\\
w_{(1)}\otimes w_{(2)}\otimes F&\mapsto& F(w_{(1)}\otimes w_{(2)})
\end{eqnarray*}
gives a natural map
$${\cal F}[Q(z)]^{W_{3}}_{W_{1}W_{2}}: W_{1}\otimes W_{2}\to
\mbox{\rm Hom}({\cal M}[Q(z)]_{W_{1}W_{2}}^{W_{3}}, \overline{W}_{3})=
({\cal M}[Q(z)]^{W_{3}}_{W_{1}W_{2}})^{*}\otimes \overline{W}_{3}.$$
Also,
$({\cal M}[Q(z)]^{W_{3}}_{W_{1}W_{2}})^{*}\otimes W_{3}$ is a $V$-module
(with finite-dimensional weight spaces) in the obvious way, and the map
${\cal F}[Q(z)]^{W_{3}}_{W_{1}W_{2}}$
is clearly a $Q(z)$-intertwining map, where we
make the identification
$$({\cal M}[Q(z)]^{W_{3}}_{W_{1}W_{2}})^{*}\otimes \overline{W}_{3}
=\overline{({\cal M}[Q(z)]^{W_{3}}_{W_{1}W_{2}})^{*}\otimes W_{3}}.$$
This gives us a natural $Q(z)$-product.

 Now we consider a special but important class of
vertex operator algebras satisfying certain finiteness and
semisimplicity conditions.
\begin{defi}
{\rm A  vertex operator algebra $V$ is {\it rational} if it
satisfies the following conditions:
\begin{enumerate}
\item There are only finitely many irreducible $V$-modules (up to equivalence).
\item Every $V$-module is completely reducible (and is in particular a
{\it finite} direct sum of irreducible modules).
\item All the fusion rules for $V$ are finite (for triples of
irreducible modules and hence arbitrary modules).
\end{enumerate}
}
\end{defi}

The next result shows that $Q(z)$-tensor products exist
for the category of modules
for a rational vertex operator algebra.

\begin{propo}
Let $V$ be rational and let $W_{1}$, $W_{2}$ be $V$-modules. Then
$(W_{1}\boxtimes_{Q(z)}W_{2}, Y_{Q(z)}; \boxtimes_{Q(z)})$ exists, and in fact
$$W_{1}\boxtimes_{Q(z)}W_{2}=\coprod_{i=1}^{k}
({\cal M}[Q(z)]^{M_{i}}_{W_{1}W_{2}})^{*}\otimes M_{i},$$
where $\{ M_{1}, \dots, M_{k}\}$ is a set of representatives of the
equivalence classes of irreducible $V$-modules,  and the right-hand side
 is equipped with the $V$-module and $Q(z)$-product structure indicated
above. That is,
$$\boxtimes_{Q(z)}=\sum_{i=1}^{k}{\cal F}[Q(z)]^{M_{i}}_{W_{1}W_{2}}.$$
\end{propo}

This construction of a $Q(z)$-tensor product module is essentially
tautological. We now describe two constructions, which are more
useful, of a
$Q(z)$-tensor product of two modules for a vertex operator algebra $V$, in the
presence of a certain hypothesis which holds in case $V$ is rational.
Fix  a nonzero complex number $z$ and $V$-modules $(W_{1}, Y_{1})$ and
$(W_{2}, Y_{2})$ as before.
We first define a linear action of $V \otimes \iota_{+}{\Bbb C}[t,t^{-
1},(z+t)^{-1}]$ on $(W_1 \otimes W_2)^*$,  that is, a linear map
$$\tau_{Q(z)}: V\otimes \iota_{+}{\Bbb C}[t, t^{-1}, (z+t)^{-1}]\to
\mbox{\rm End}\;(W_{1}\otimes W_{2})^{*},$$
where $\iota_+$ is the operation of expanding a rational function in
the formal variable $t$ in the direction of positive powers of $t$
(with at most finitely many negative powers of $t$), by
\begin{eqnarray*}
\lefteqn{\left(\tau_{Q(z)}
\left(z^{-1}\delta\left(\frac{x_{1}-x_{0}}{z}\right)
Y_{t}(v, x_{0})\right)\lambda\right)(w_{(1)}\otimes w_{(2)})}\nonumber\\
&&=x^{-1}_{0}\delta\left(\frac{x_{1}-z}{x_{0}}\right)
\lambda(Y_{1}^{*}(v, x_{1})w_{(1)}\otimes w_{(2)})\nonumber\\
&&\hspace{2em}-x_{0}^{-1}\delta\left(\frac{z-x_{1}}{-x_{0}}\right)
\lambda(w_{(1)}\otimes Y_{2}(v, x_{1})w_{(2)}).
\end{eqnarray*}
(Recall that $Y_t$ and $Y^*$ have been defined above, and note that
the coefficients of the monomials in $x_0$ and $x_1$ in
$$z^{-1}\delta\left(\frac{x_{1}-x_{0}}{z}\right)Y_{t}(v,x_0),$$
for all $v \in V,$ span the space
$V\otimes \iota_{+}{\Bbb C}[t, t^{-1}, (z+t)^{-1}]$.)
Write
$$Y'_{Q(z)}(v, x)=\tau_{Q(z)}(Y_{t}(v, x)).$$

\begin{propo}
The action $Y'_{Q(z)}$ satisfies the commutator formula for vertex operators,
that is, on
$(W_{1}\otimes W_{2})^{*}$,
\begin{eqnarray*}
\lefteqn{[Y'_{Q(z)}(v_{1}, x_{1}), Y'_{Q(z)}(v_{2}, x_{2})]=}\nno\\
&&=\res_{x_{0}}x_{2}^{-1}\delta\left(\frac{x_{1}-x_{0}}{x_{2}}\right)
Y'_{Q(z)}(Y(v_{1}, x_{0})v_{2}, x_{2})
\end{eqnarray*}
for $v_{1}, v_{2}\in V$.
\end{propo}

(The symbol $\res_x$ designates the coefficient of $x^{-1}$ in a
formal series.)

Let  $W_{3}$ be another $V$-module.
Note that $V\otimes \iota_{+}{\Bbb C}[t, t^{-1}, (z+t)^{-1}]$ acts on
$W'_{3}$ in the obvious way. The following  result
provides some important motivation for the definition of our action on
$(W_{1}\otimes W_{2})^{*}$:

\begin{propo}\label{5-3}
Under the natural isomorphism
\begin{equation}\label{5.15}
\mbox{\rm Hom}(W'_{3}, (W_{1}\otimes
W_{2})^{*})\stackrel{\sim}{\to}\mbox{\rm Hom}(W_{1}\otimes W_{2},
\overline{W}_{3}),
\end{equation}
the maps in $\mbox{\rm Hom}(W'_{3}, (W_{1}\otimes
W_{2})^{*})$ intertwining the two actions of
$V \otimes \iota_{+}{\Bbb C}[t,t^{-
1},(z+t)^{-1}]$ on $W'_{3}$ and
$(W_{1}\otimes W_{2})^{*}$ correspond exactly to the
$Q(z)$-intertwining maps of type ${W_{3}}\choose {W_{1}W_{2}}$.
\end{propo}

\begin{rema}
{\rm Combining the last result with Proposition \ref{4-7}, we see that
the maps in $\mbox{\rm Hom}(W'_{3}, (W_{1}\otimes
W_{2})^{*})$ intertwining the two actions on $W'_{3}$ and
$(W_{1}\otimes W_{2})^{*}$ also correspond exactly to the
intertwining operators of type
${W_{1}'}\choose {W'_{3}\;W_{2}}$. In particular, given any integer $p$,
the map
$(F^{Q(z)}_{{\cal Y}, p})': W'_{3}\to (W_{1}\otimes W_{2})^{*}$ defined by
$$
(F^{Q(z)}_{{\cal Y}, p})'(w'_{(3)})(w_{(1)}\otimes w_{(2)})=\bra w_{(1)},
{\cal Y}(w'_{(3)}, e^{l_{p}(z)})w_{(2)}\ket_{W'_1}
$$
(recall (\ref{4.9}))  intertwines the
action of $V\otimes \iota_{+}{\Bbb C}[t, t^{-1}, (z+t)^{-1}]$ on $W'_{3}$ and
the action $\tau_{Q(z)}$ of $V\otimes \iota_{+}{\Bbb C}[t, t^{-1}, (z+t)^{-1}]$
on $(W_{1}\otimes W_{2})^{*}$.}
\end{rema}

Consider the following nontrivial and subtle {\it compatibility
condition}
on $\lambda
\in (W_{1}\otimes W_{2})^{*}$: The formal Laurent series $Y'_{Q(z)}(v,
x_{0})\lambda$ involves only finitely many negative powers of $x_{0}$
and
\begin{eqnarray}\label{5.18}
\lefteqn{\tau_{Q(z)}\left(z^{-1}\delta\left(\frac{x_{1}-x_{0}}{z}\right)
Y_{t}(v, x_{0})\right)\lambda=}\nno\\
&&=z^{-1}\delta\left(\frac{x_{1}-x_{0}}{z}\right)
Y'_{Q(z)}(v, x_{0})\lambda  \;\;\;\;\; \mbox{\rm for all}\;\;v\in V.
\end{eqnarray}
(Note that the two sides are not {\it a priori} equal for general
$\lambda\in (W_{1}\otimes W_{2})^{*}$.)

Let $W$ be a subspace of $(W_{1}\otimes W_{2})^{*}$.  We say that $W$
is {\it compatible for $\tau_{Q(z)}$} if every element of $W$
satisfies the compatibility condition. Also, we say that $W$ is
(${\Bbb C}$-){\it graded} if it is ${\Bbb C}$-graded by its weight
subspaces, and that $W$ is a $V$-{\it module} (respectively, {\it
generalized module}) if $W$ is graded and is a module (respectively,
generalized module) when equipped with this grading and with the
action of $Y'_{Q(z)}(\cdot, x)$.  (The notion of ``weight'' for
$(W_{1}\otimes W_{2})^{*}$ is defined by means of the eigenvalues of
the operator $L'_{Q(z)}(0) = \res_x Y'_{Q(z)}(\omega,x)$; recall that
$\omega$ is the element of $V$ giving the Virasoro algebra.)  A sum of
compatible modules or generalized modules is clearly a generalized
module. The weight subspace of a subspace $W$ with weight $n\in {\Bbb
C}$ will be denoted $W_{(n)}$.

Define
$$W_{1}\hboxtr_{Q(z)}W_{2}=\sum_{W\in {\cal W}_{Q(z)}}W =\bigcup_{W\in
{\cal W}_{Q(z)}} W\subset
(W_{1}\otimes W_{2})^{*},$$
where ${\cal W}_{Q(z)}$ is the set all compatible
modules for $\tau_{Q(z)}$ in $(W_{1}\otimes W_{2})^{*}$.
Then $W_{1}\hboxtr_{Q(z)}W_{2}$
is a compatible generalized module. We have:

\begin{propo}\label{5-5}
The subspace $W_{1}\hboxtr_{Q(z)}W_{2}$ of $(W_{1}\otimes W_{2})^{*}$ is
a generalized module with the following
property: Given  any $V$-module $W_{3}$,
there is a natural linear isomorphism determined by (\ref{5.15})
between the space of
all $Q(z)$-intertwining maps of type ${W_{3}}\choose {W_{1}\;W_{2}}$
and the space of all maps of generalized modules {from} $W'_{3}$ to
$W_{1}\hboxtr _{Q(z)}W_{2}$.
\end{propo}

\begin{propo}\label{5-6}
Let $V$ be a rational vertex operator algebra and $W_{1}$, $W_{2}$ two
$V$-modules. Then $W_{1}\hboxtr _{Q(z)}W_{2}$ is a module.
\end{propo}

Now we assume that $W_{1}\hboxtr_{Q(z)} W_{2}$ is a module (which
occurs if $V$ is rational, by
the last proposition).  In this case, we define a $V$-module
$W_{1}\boxtimes_{Q(z)} W_{2}$ by
\begin{equation}\label{5.20}
W_{1}\boxtimes_{Q(z)} W_{2}=(W_{1}\hboxtr_{Q(z)}W_{2})'
\end{equation}
(note that by our choice of notation, $\hboxtr'=\boxtimes$)
and we write the corresponding
action as $Y_{Q(z)}$.
Applying Proposition \ref{5-5} to the special module
$W_{3}=W_{1}\boxtimes_{Q(z)} W_{2}$ and the natural isomorphism {}from
the double contragredient module $W'_{3}=(W_{1}\hboxtr_{Q(z)}
W_{2})''$
to
$W_{1}\hboxtr_{Q(z)} W_{2}$ (recall \cite{FHL}, Theorem 5.3.1),
we obtain using
(\ref{5.15}) a canonical $Q(z)$-intertwining map of type
${W_{1}\boxtimes_{Q(z)} W_{2}}\choose {W_{1}W_{2}}$, which we denote
\begin{eqnarray*}
\boxtimes_{Q(z)}: W_{1}\otimes W_{2}&\to &
\overline{W_{1}\boxtimes_{Q(z)} W_{2}}\nno\\
w_{(1)}\otimes  w_{(2)}&\mapsto& w_{(1)} \boxtimes_{Q(z)}w_{(2)}.
\end{eqnarray*}
This is the unique linear map such that
$$\langle \lambda, w_{(1)} \boxtimes_{Q(z)}
w_{(2)}\rangle_{W_{1}\boxtimes_{Q(z)} W_{2}}
=\lambda(w_{(1)}\otimes w_{(2)})$$
for all $w_{(1)}\in W_{1}$, $w_{(2)}\in W_{2}$ and
$\lambda\in W_{1}\hboxtr_{Q(z)} W_{2}$. Moreover, we have:

\begin{propo}\label{5-7}
The $Q(z)$-product $(W_{1}\boxtimes_{Q(z)} W_{2}, Y_{Q(z)};
\boxtimes_{Q(z)})$
is a $Q(z)$-tensor product of $W_{1}$ and $W_{2}$.
\end{propo}

More generally, dropping the assumption that $W_{1}\hboxtr_{Q(z)}
W_{2}$ is a module, we have:

\begin{propo}\label{5-8}
The $Q(z)$-tensor product of $W_{1}$ and $W_{2}$ exists (and is given
by (\ref{5.20})) if and only if $W_{1}\hboxtr_{Q(z)} W_{2}$ is a module.
\end{propo}

It is not difficult to see that any element of
$W_{1}\hboxtr_{Q(z)} W_{2}$ is an element $\lambda$
of $(W_{1}\otimes W_{2})^{*}$ satisfying:

\begin{description}
\item[The compatibility condition] (recall (\ref{5.18})){\bf :}
{\bf (a)} The  {\it lower
truncation condition}:
For all $v\in V$, the formal Laurent series $Y'_{Q(z)}(v, x)\lambda$
involves only finitely many negative
powers of $x$.

{\bf (b)} The following formula holds:
\begin{eqnarray*}
\lefteqn{\tau_{Q(z)}\left(z^{-1}\delta\left(\frac{x_{1}-x_{0}}{z}\right)
Y_{t}(v, x_{0})\right)\lambda=}\nno\\
&&=z^{-1}\delta\left(\frac{x_{1}-x_{0}}{z}\right)
Y'_{Q(z)}(v, x_{0})\lambda  \;\;\;\;\; \mbox{\rm for all}\;\;v\in V.
\end{eqnarray*}

\item[The local grading-restriction  condition:]
{\bf (a)} The {\it grading condition}:
$\lambda$ is a (finite) sum of
weight vectors of $(W_{1}\otimes W_{2})^{*}$.

{\bf (b)} Let $W_{\lambda}$ be the smallest subspace of $(W_{1}\otimes
W_{2})^{*}$ containing $\lambda$ and stable under the component
operators $\tau_{Q(z)}(v\otimes t^{n})$ of the operators $Y'_{Q(z)}(v,
x)$ for $v\in V$, $n\in {\Bbb Z}$. Then the weight spaces
$(W_{\lambda})_{(n)}$, $n\in {\Bbb C}$, of the (graded) space
$W_{\lambda}$ have the properties
\begin{eqnarray*}
&\mbox{\rm dim}\ (W_{\lambda})_{(n)}<\infty \;\;\;\mbox{\rm for}\
n\in {\Bbb C},&\\
&(W_{\lambda})_{(n)}=0 \;\;\;\mbox{\rm for $n$ whose real part is
sufficiently small.}&
\end{eqnarray*}
\end{description}

We have the following basic ``second construction'' of
$W_{1}\hboxtr_{Q(z)} W_{2}$:

\begin{theo}\label{6-3}
The subspace of $(W_{1}\otimes W_{2})^{*}$ consisting of the
elements satisfying the compatibility
condition and the local grading-restriction condition, equipped with
$Y'_{Q(z)}$, is a generalized module and is equal to
$W_{1}\hboxtr_{Q(z)} W_{2}$.
\end{theo}

The following result follows immediately {from} Proposition \ref{5-8}, the
theorem above and the definition of $W_{1}\boxtimes_{Q(z)} W_{2}$:

\begin{corol}\label{6-4}
The $Q(z)$-tensor product of $W_{1}$ and $W_{2}$ exists if and only if
the subspace of $(W_{1}\otimes W_{2})^{*}$ consisting of the elements
satisfying the compatibility condition and the local
grading-restriction condition, equipped with $Y'_{Q(z)}$,
is a module. In this case, this module coincides with the module $W_1
\hboxtr_{Q(z)} W_2$, and the contragredient module of this
module, equipped with the $Q(z)$-intertwining map $\boxtimes_{Q(z)}$,
 is a $Q(z)$-tensor product of $W_{1}$ and
$W_{2}$,  equal to the structure $(W_{1}\boxtimes_{Q(z)}
W_{2}, Y_{Q(z)};
\boxtimes_{Q(z)})$ constructed above.
\end{corol}

{From} this result and Propositions \ref{5-6} and \ref{5-7}, we have:

\begin{corol}
Let $V$ be a rational vertex operator algebra and $W_{1}$, $W_{2}$ two
$V$-modules. Then the $Q(z)$-tensor product
$(W_{1}\boxtimes_{Q(z)} W_{2}, Y_{Q(z)};
\boxtimes_{Q(z)})$ may be constructed as described in Corollary \ref{6-4}.
\end{corol}

This finishes our review of the results on $Q(z)$-tensor products,
{}from Parts I and II.

The numberings of sections, formulas, etc., in Part III
continue those of Parts I and II.

Part III contains two sections following this introduction.  In
Section 12, we recall the notion of $P(z)$-tensor product introduced
in Part I, and we present the basic properties of this tensor product
analogous to those of the $Q(z)$-tensor product, including the
existence of $P(z)$-tensor products in the case of rational vertex
operator algebras.  We only state the results since all the proofs are
exactly analogous to those for $Q(z)$-tensor products (see Parts I and
II).  We establish two constructions of a $P(z)$-tensor product,
analogous to those of the $Q(z)$-tensor product, in Section 13.

\paragraph{Acknowledgments}
During the
course of this work, Y.-Z.~H. has been supported in part by NSF grants
DMS-9104519 and DMS-9301020 and J.~L. by NSF grant
DMS-9111945 and the Rutgers University Faculty
Academic Study Program.  J.~L. also thanks the Institute for Advanced
Study for its hospitality in fall, 1992.

\setcounter{section}{11}
\renewcommand{\theequation}{\thesection.\arabic{equation}}
\renewcommand{\therema}{\thesection.\arabic{rema}}
\setcounter{equation}{0}
\setcounter{rema}{0}

\section{Properties of $P(z)$-tensor products}

We first review the notion of $P(z)$-tensor product introduced in \cite{HL1}.
Fix $z\in {\Bbb C}^{\times}$.  Recall that $V$ is a fixed vertex
operator algebra and let $(W_{1}, Y_{1})$, $(W_{2}, Y_{2})$ and
 $(W_{3}, Y_{3})$ be $V$-modules.
  By a {\it $P(z)$-intertwining map of type ${W_{3}}\choose
{W_{1}W_{2}}$} we mean  a linear map $F: W_{1}\otimes W_{2} \to
\overline{W}_{3}$ satisfying the condition
\begin{eqnarray}
\lefteqn{x_{0}^{-1}\delta\left(\frac{ x_{1}-z}{x_{0}}\right)
Y_{3}(v, x_{1})F(w_{(1)}\otimes w_{(2)})=}\nonumber\\
&&=z^{-1}\delta\left(\frac{x_{1}-x_{0}}{z}\right)
F(Y_{1}(v, x_{0})w_{(1)}\otimes w_{(2)})\nonumber\\
&&\hspace{2em}+x_{0}^{-1}\delta\left(\frac{z-x_{1}}{-x_{0}}\right)
F(w_{(1)}\otimes Y_{2}(v, x_{1})w_{(2)})
\end{eqnarray}
for $v\in V$, $w_{(1)}\in W_{1}$, $w_{(2)}\in W_{2}$.
(The expressions in (12.1) are well defined, as explained in Part I.)

We denote the vector space of $P(z)$-intertwining maps of type
${W_{3}}\choose {W_{1}W_{2}}$ by ${\cal M}[P(z)]^{W_{3}}_{W_{1}W_{2}}$.

A {\it $P(z)$-product of $W_{1}$
and $W_{2}$} is a $V$-module $(W_{3}, Y_{3})$ equipped with a
$P(z)$-intertwining map $F$ of type ${W_{3}}\choose {W_{1}W_{2}}$.  We
denote it by $(W_{3}, Y_{3}; F)$ (or simply by $(W_{3}, F)$).  Let
$(W_{4}, Y_{4}; G)$ be another $P(z)$-product of $W_{1}$ and $W_{2}$.
A {\it morphism} {}from $(W_{3}, Y_{3}; F)$ to $(W_{4}, Y_{4}; G)$ is a
module map $\eta$ {}from $W_{3}$ to $W_{4}$ such that
\begin{equation}
G=\overline{\eta}\circ F,
\end{equation}
where $\overline{\eta}$ is the natural map {}from $\overline{W}_{3}$ to
$\overline{W}_{4}$ uniquely extending $\eta$.

\begin{defi}
{\rm A {\it $P(z)$-tensor product of $W_{1}$ and $W_{2}$} is
a $P(z)$-product $$(W_{1}\boxtimes_{P(z)} W_{2}, Y_{P(z)};
\boxtimes_{P(z)})$$
such that for any $P(z)$-product
$(W_{3}, Y_{3}; F)$, there is a unique morphism {}from
$$(W_{1}\boxtimes_{P(z)} W_{2}, Y_{P(z)};
\boxtimes_{P(z)})$$ to $(W_{3}, Y_{3}; F)$.
The $V$-module $(W_{1}\boxtimes_{P(z)} W_{2},  Y_{P(z)})$ is
called a {\it $P(z)$-tensor product module} of $W_{1}$ and $W_{2}$.}
\end{defi}

It is clear that a $P(z)$-tensor product of $W_1$ and $W_2$ is unique
up to unique isomorphism if it exists.

We now describe the precise connection between
 intertwining operators
and $P(z)$-intertwining maps of the same type. Fix an integer $p$.
Let ${\cal Y}$ be an intertwining operator of
type ${W_{3}}\choose {W_{1}W_{2}}$. We have a linear map
$F^{P(z)}_{{\cal Y}, p}: W_{1}\otimes W_{2}\to \overline{W}_{3}$ given by
\begin{equation}
F^{P(z)}_{{\cal Y}, p}(w_{(1)}\otimes w_{(2)})
={\cal Y}(w_{(1)}, e^{l_{p}(z)})w_{(2)}
\end{equation}
for all $w_{(1)}\in W_{1}$, $w_{(2)}\in W_{2}$. Using the Jacobi identity for
${\cal Y}$, we see easily that $F^{P(z)}_{{\cal Y}, p}$
is a $P(z)$-intertwining map.
Conversely, given  a $P(z)$-intertwining map $F$, homogeneous elements
$w_{(1)}\in W_{1}$ and $w_{(2)}\in W_{2}$ and $n\in {\Bbb C}$, we define
$(w_{(1)})_{n}w_{(2)}$ to be the projection of the image of $w_{(1)}\otimes
w_{(2)}$ under $F$ to the homogeneous subspace of $W_{3}$ of weight
$$\mbox{\rm wt}\ w_{(1)} -n-1+\mbox{\rm wt}\ w_{(2)},$$
multiplied by $e^{(n+1)l_p(z)}$.
Using this, we define
\begin{equation}
{\cal Y}_{F, p}(w_{(1)}, x)w_{(2)}=\sum_{n\in
{\Bbb C}}(w_{(1)})_{n}w_{(2)}x^{-n-1},
\end{equation}
and by linearity, we obtain a linear map
\begin{eqnarray*}
W_{1}\otimes W_{2}&\to& W_{3}\{ x\}\\
w_{(1)}\otimes w_{(2)}&\mapsto &{\cal Y}_{F, p}(w_{(1)}, x)w_{(2)}.
\end{eqnarray*}

The proof of the following result is analogous to (and slightly shorter than)
the proof of the corresponding result for $Q(z)$-intertwining maps in
\cite{HL1}:

\begin{propo}
For $p\in {\Bbb Z}$, the correspondence ${\cal Y}\mapsto
F^{P(z)}_{{\cal Y},p}$ is a
linear isomorphism
{from} the vector space ${\cal V}^{W_{3}}_{W_{1}W_{2}}$ of
intertwining operators of type ${W_{3}}\choose {W_{1}\; W_{2}}$
to the vector space ${\cal M}[P(z)]^{W_{3}}_{W_{1}W_{2}}$ of
$P(z)$-intertwining maps of
type ${W_{3}}\choose {W_{1}W_{2}}$. Its inverse map is given by
$F\mapsto {\cal Y}_{F, p}$.
\epf
\end{propo}

The following immediate result relates module maps {from} a $P(z)$-tensor
product module with intertwining maps and intertwining operators:
\begin{propo}
Suppose that $W_{1}\boxtimes_{P(z)}W_{2}$ exists. We have a natural
isomorphism
\begin{eqnarray}
\mbox{\rm Hom}_{V}(W_{1}\boxtimes_{P(z)}W_{2}, W_{3})&
\stackrel{\sim}{\to}&
{\cal M}[P(z)]^{W_{3}}_{W_{1}W_{2}}\nno\\
\eta&\mapsto& \overline{\eta}\circ \boxtimes_{P(z)}
\end{eqnarray}
 and for $p\in {\Bbb Z}$, a
natural isomorphism
\begin{eqnarray}
\mbox{\rm Hom}_{V}(W_{1}\boxtimes_{P(z)} W_{2}, W_{3})&
\stackrel{\sim}{\rightarrow}&
{\cal V}^{W_{3}}_{W_{1}W_{2}}\nno\\
\eta&\mapsto & {\cal Y}_{\eta, p}
\end{eqnarray}
where ${\cal Y}_{\eta, p}={\cal Y}_{F, p}$ with
$F=\overline{\eta}\circ \boxtimes_{P(z)}$.
\epf
\end{propo}

It is clear {from} the Definition 12.1 that $P(z)$-tensor product
 that the $P(z)$-tensor product operation
distributes over direct sums in the following sense:

\begin{propo}
Let $U_{1}, \dots, U_{k}$, $W_{1}, \dots, W_{l}$ be $V$-modules and suppose
that each $U_{i}\boxtimes_{P(z)}W_{j}$ exists. Then
$(\coprod_{i}U_{i})\boxtimes_{P(z)}(\coprod_{j}W_{j})$ exists and there is a
natural isomorphism
\begin{equation}
\biggl(\coprod_{i}U_{i}\biggr)\boxtimes_{P(z)}\biggl(\coprod_{j}W_{j}\biggr)
\stackrel{\sim}
{\rightarrow} \coprod_{i,j}U_{i}\boxtimes_{P(z)}W_{j}.
\hspace{2em}\Box
\end{equation}
\end{propo}

Now consider $V$-modules $W_{1}$, $W_{2}$ and $W_{3}$ and
suppose that
$$\mbox{\rm dim}\ {\cal V}^{W_{3}}_{W_{1}W_{2}} <\infty \ \;
(\mbox{\rm or}\;\ \mbox{\rm dim}\ {\cal M}[P(z)]^
{W_{3}}_{W_{1}W_{2}} <\infty).$$
The natural evaluation map
\begin{eqnarray}
W_{1}\otimes W_{2}\otimes {\cal M}[P(z)]^{W_{3}}_{W_{1}W_{2}}&\to&
\overline{W}_{3}
\nno\\
w_{(1)}\otimes w_{(2)}\otimes F&\mapsto& F(w_{(1)}\otimes w_{(2)})
\end{eqnarray}
gives a natural map
\begin{equation}
{\cal F}[P(z)]^{W_{3}}_{W_{1}W_{2}}: W_{1}\otimes W_{2}\to
({\cal M}[P(z)]^{W_{3}}_{W_{1}W_{2}})^{*}\otimes \overline{W}_{3}.
\end{equation}
Since dim$\ {\cal M}[P(z)]^{W_{3}}_{W_{1}W_{2}}<\infty$,
$({\cal M}[P(z)]^{W_{3}}_{W_{1}W_{2}})^{*}\otimes W_{3}$ is a $V$-module
(with finite-dimensional weight spaces) in the obvious way, and the map
${\cal F}[P(z)]^{W_{3}}_{W_{1}W_{2}}$ is clearly a $P(z)$-intertwining map,
where we
make the identification
\begin{equation}
({\cal M}[P(z)]^{W_{3}}_{W_{1}W_{2}})^{*}\otimes \overline{W}_{3}
=\overline{({\cal M}[P(z)]^{W_{3}}_{W_{1}W_{2}})^{*}\otimes W_{3}}.
\end{equation}
This gives us a natural $P(z)$-product.

The next result  shows that $P(z)$-tensor products exist for the category of
modules
for a rational vertex operator algebra (recall Definition I.6). It is
proved by the same argument used to prove the analogous result
in \cite{HL1} for the $Q(z)$-tensor product. As in the case of $Q(z)$-tensor
products, there is no need to assume that
$W_{1}$ and $W_{2}$ are irreducible in the formulation or proof,
 but by Proposition 12.4, the case in which $W_{1}$ and $W_{2}$ are irreducible
gives all the necessary information, and the tensor product is canonically
described using only the spaces of intertwining maps among triples of
{\it irreducible} modules.

\begin{propo}
Let $V$ be rational and let $W_{1}$, $W_{2}$ be $V$-modules. Then
$$(W_{1}\boxtimes_{P(z)}W_{2}, Y_{P(z)}; \boxtimes_{P(z)})$$ exists, and in
fact
\begin{equation}
W_{1}\boxtimes_{P(z)}W_{2}=\coprod_{i=1}^{k}
({\cal M}[P(z)]^{M_{i}}_{W_{1}W_{2}})^{*}\otimes M_{i},
\end{equation}
where $\{ M_{1}, \dots, M_{k}\}$ is a set of representatives of the
equivalence classes of irreducible $V$-modules,  and the right-hand side of
(12.11) is equipped with the $V$-module and $P(z)$-product structure indicated
above. That is,
\begin{equation}
\boxtimes_{P(z)}=\sum_{i=1}^{k}{\cal F}[P(z)]^{M_{i}}_{W_{1}W_{2}}.
\hspace{2em}\Box
\end{equation}
\end{propo}

\begin{rema}
By combining Proposition 12.5 with Proposition 12.2, we can
express $W_{1}\boxtimes_{P(z)} W_{2}$ in terms of
${\cal V}^{M_{i}}_{W_{1}W_{2}}$  in
place of ${\cal M}[P(z)]^{M_{i}}_{W_{1}W_{2}}$.
\end{rema}

The construction in Proposition 12.5 is tautological, and we view the
argument as essentially an existence proof. We shall give two
constructions of a $P(z)$-tensor product, under suitable conditions,
in the next section.

\renewcommand{\theequation}{\thesection.\arabic{equation}}
\renewcommand{\therema}{\thesection.\arabic{rema}}
\setcounter{equation}{0}
\setcounter{rema}{0}

\section{Constructions of $P(z)$-tensor product}

In this section, using the results of \cite{HL1} and \cite{HL2}
(see also the introduction to the present paper),
we give two constructions of a $P(z)$-tensor product of two
$V$-modules $W_{1}$ and $W_{2}$ when certain conditions are satisfied.
This treatment is parallel to that of
the $Q(z)$-tensor product $W_{1}\boxtimes_{Q(z)} W_{2}$
of $W_{1}$ and $W_{2}$. In particular,
when the vertex operator algebra is rational, we construct a $P(z)$-tensor
product of $W_{1}$ and $W_{2}$ in ways that are more useful than
Proposition 12.5.

Combining Proposition \ref{4-7} (Proposition 4.7 in \cite{HL1}),
Proposition \ref{4-9} (Proposition 4.9 in \cite{HL1}) and Proposition
12.2, we obtain an isomorphism {}from the space ${\cal
M}[Q(z)]_{W_{1}W_{2}}^{W_{1}\boxtimes_{Q(z)} W_{2}}$ of
$Q(z)$-intertwining maps to the space ${\cal
M}[P(z)]_{W_{1}W_{2}}^{W_{1}\boxtimes_{Q(z)} W_{2}}$ of
$P(z)$-intertwining maps for any pair $(p_{1}, p_{2})$ of integers,
when $W_{1}\boxtimes_{Q(z)}W_{2}$ exists.  Thus for any pair $(p_{1},
p_{2})$ of integers, we have a $P(z)$-product consisting of the module
$W_{1}\boxtimes_{Q(z)} W_{2}$ and the $P(z)$-intertwining map which is
the image of the $Q(z)$-intertwining map $\boxtimes_{Q(z)}$ under the
isomorphism corresponding to $(p_{1}, p_{2})$.  It is easy to show
that this is a $P(z)$-tensor product of $W_{1}$ and $W_{2}$.  Since we
shall be interested in associativity and other nice properties, we
shall not discuss in detail the proof that the $P(z)$-intertwining map
above gives a $P(z)$-tensor product. Instead, we would like to
construct a $P(z)$-tensor product in a way analogous to the
construction of the $Q(z)$-tensor product in Section 5 and 6. To
construct such a $P(z)$-tensor product of two modules, we could first
prove results analogous to those used in the constructions of the
$Q(z)$-tensor product, and then obtain the $P(z)$-tensor product. But
there is an easier way: We can use some of the results proved for the
construction of the $Q(z^{-1})$-tensor product (mainly the results in
Section 9, 10 and 11) to derive the results that we want for the
construction of the $P(z)$-tensor product. This is what we shall do in
this section.  We chose to construct the $Q(z)$-tensor product first
in the present series of papers because in the case that $z=1$ and
$W_{1}$, $W_{2}$ are modules for a vertex operator algebra associated
to an affine Lie algebra, it can be proved that the $Q(1)$-tensor
product of $W_{1}$ and $W_{2}$ agrees with the tensor product
constructed by Kazhdan and Lusztig in \cite{KL}. However, {}from the
geometric viewpoint, the simplest case of the associativity is the one
for $P(z)$-tensor products ($z\in {\Bbb C}^{\times}$), although in the
proof of this associativity, $Q(z)$-tensor products are also used
\cite{H}. The interested reader can imitate the constructions and
calculations in \cite{HL1} and
\cite{HL2} to obtain the $P(z)$-tensor product directly. We emphasize,
though, that the $P(z)$- and $Q(z)$-tensor products are on equal
footing.

For two $V$-modules $(W_{1}, Y_{1})$ and $(W_{2}, Y_{2})$,
we define a linear action of $$V \otimes \iota_{+}{\Bbb C}[t,t^{- 1},
(z^{-1}-t)^{-1}]$$ on $(W_1 \otimes W_2)^*$ (recall that
 $\iota_{+}$ denotes the operation of expansion of a rational
function of $t$ in the direction of positive powers of $t$),  that is, a
linear map
\begin{equation}
\tau_{P(z)}: V\otimes \iota_{+}{\Bbb C}[t, t^{-1}, (z^{-1}-t)^{-1}]\to
\mbox{\rm End}\;(W_{1}\otimes W_{2})^{*},
\end{equation}
by
\begin{eqnarray}
\lefteqn{\left(\tau_{P(z)}
\left(x_{0}^{-1}\delta\left(\frac{x^{-1}_{1}-z}{x_{0}}\right)
Y_{t}(v, x_{1})\right)\lambda\right)(w_{(1)}\otimes w_{(2)})=}\nonumber\\
&&=z^{-1}\delta\left(\frac{x^{-1}_{1}-x_{0}}{z}\right)
\lambda(Y_{1}(e^{x_{1}L(1)}(-x_{1}^{-2})^{L(0)}v, x_{0})w_{(1)}\otimes w_{(2)})
\nonumber\\
&&\quad +x^{-1}_{0}\delta\left(\frac{z-x^{-1}_{1}}{-x_{0}}\right)
\lambda(w_{(1)}\otimes Y_{2}^{*}(v, x_{1})w_{(2)})
\end{eqnarray}
for $v\in V$, $\lambda\in (W_{1}\otimes W_{2})^{*}$, $w_{(1)}\in W_{1}$,
$w_{(2)}\in W_{2}$.
The formula (13.2) does indeed give a well-defined map of the type
(13.1) (in generating-function form); this definition is motivated by
(12.1), in which we first replace $x_1$ by $x_1^{-1}$ and then replace
$v$ by $e^{x_{1}L(1)}(-x_{1}^{-2})^{L(0)}v$.

Let
$W_{3}$ be another $V$-module.
The space $V\otimes \iota_{+}{\Bbb C}[t, t^{-1}, (z^{-1}-t)^{-1}]$ acts on
$W'_{3}$ in the obvious way, where $v\otimes t^{n}$ ($v\in V$, $n\in {\Bbb Z}$)
acts as the component $v_{n}$ of $Y(v, x)$. The following  result, which
follows
immediately {from} the definitions (12.1) and (13.2) and is analogous to
Proposition \ref{5-3} (Proposition 5.3 in \cite{HL1}),
provides further motivation for the definition of our action on
$(W_{1}\otimes W_{2})^{*}$:

\begin{propo}
Under the natural isomorphism
\begin{equation}
\mbox{\rm Hom}(W'_{3}, (W_{1}\otimes
W_{2})^{*})\stackrel{\sim}{\to}\mbox{\rm Hom}(W_{1}\otimes W_{2},
\overline{W}_{3}),
\end{equation}
the maps in $\mbox{\rm Hom}(W'_{3}, (W_{1}\otimes
W_{2})^{*})$ intertwining the two actions of
$$V \otimes \iota_{+}{\Bbb C}[t,t^{-
1},(z^{-1}-t)^{-1}]$$ on $W'_{3}$ and
$(W_{1}\otimes W_{2})^{*}$ correspond exactly to the
$P(z)$-intertwining maps of type ${W_{3}}\choose {W_{1}W_{2}}$. \epf
\end{propo}

\begin{rema}
{\rm Combining the last result with Proposition 12.2, we see that
the maps in $\mbox{\rm Hom}(W'_{3}, (W_{1}\otimes
W_{2})^{*})$ intertwining the two actions on $W'_{3}$ and
$(W_{1}\otimes W_{2})^{*}$ also correspond exactly to the
intertwining operators of type
${W_{3}}\choose {W_{1}\;W_{2}}$. In particular, for any intertwining operator
${\cal Y}$ of type ${W_{3}}\choose {W_{1}W_{2}}$ and any integer $p$,
the map
$F'_{{\cal Y}, p}: W'_{3}\to (W_{1}\otimes W_{2})^{*}$ defined by
\begin{equation}
F'_{{\cal Y}, p}(w'_{(3)})(w_{(1)}\otimes w_{(2)})=\bra w'_{(3)},
{\cal Y}(w_{(1)}, e^{l_{p}(z)})w_{(2)}\ket_{W_{3}}
\end{equation}
for $w'_{(3)}\in W'_{3}$ (recall (12.3))  intertwines the
actions  of
$V\otimes \iota_{+}{\Bbb C}[t, t^{-1}, (z^{-1}-t)^{-1}]$ on $W'_{3}$ and
$(W_{1}\otimes W_{2})^{*}$.}
\end{rema}

Write
\begin{equation}
Y'_{P(z)}(v, x)=\tau_{P(z)}(Y_{t}(v, x))
\end{equation}
(the specialization of (13.2) to $V\otimes {\Bbb C}[t, t^{-1}]$).
More explicitly (as in formula (5.4) of Part I), (13.2) gives:
\begin{eqnarray}
\lefteqn{(Y'_{P(z)}(v, x)\lambda)(w_{(1)}\otimes w_{(2)})=}\nonumber\\
&&=\res_{x_{0}}z^{-1}\delta\left(\frac{x^{-1}-x_{0}}{z}\right)
\lambda(Y_{1}(e^{xL(1)}(-x^{-2})^{L(0)}v, x_{0})w_{(1)}\otimes w_{(2)})
\nonumber\\
&&\quad +\res_{x_{0}}x^{-1}_{0}\delta\left(\frac{z-x^{-1}}{-x_{0}}\right)
\lambda(w_{(1)}\otimes Y_{2}^{*}(v, x)w_{(2)})\nno\\
&&=\res_{x_{0}}z^{-1}\delta\left(\frac{x^{-1}-x_{0}}{z}\right)
\lambda(Y_{1}(e^{xL(1)}(-x^{-2})^{L(0)}v, x_{0})w_{(1)}\otimes w_{(2)})
\nonumber\\
&&\quad +\lambda(w_{(1)}\otimes Y_{2}^{*}(v, x)w_{(2)})
\end{eqnarray}

We have the following straightforward result, as in Proposition 5.1 of
Part I:

\begin{propo}
The action $Y'_{P(z)}$ of $V\otimes {\Bbb C}[t, t^{-1}]$ on $(W_{1}\otimes
W_{2})^{*}$ has the property
\begin{equation}
Y'_{P(z)}({\bf 1}, x)=1,
\end{equation}
where $1$ on the right-hand side is the identity map of
$(W_{1}\otimes W_{2})^{*}$, and the $L(-1)$-derivative property
\begin{equation}
\frac{d}{dx}Y'_{P(z)}(v, x)=Y'_{P(z)}(L(-1)v, x)
\end{equation}
for $v\in V$.
\end{propo}
\pf
The first part follows directly {}from the definition. We prove the
$L(-1)$-derivative property. {}From (13.6), we obtain
\begin{eqnarray}
\lefteqn{\left(\frac{d}{dx}Y'_{P(z)}(v, x)\lambda\right)
(w_{(1)}\otimes w_{(2)})=}\nno\\
&&=\frac{d}{dx}\res_{x_{0}}z^{-1}\delta\left(\frac{x^{-1}-x_{0}}{z}\right)
\lambda(Y_{1}(e^{xL(1)}(-x^{-2})^{L(0)}v, x_{0})w_{(1)}\otimes w_{(2)})
\nonumber\\
&&\quad +\frac{d}{dx}\lambda(w_{(1)}\otimes Y_{2}^{*}(v, x)w_{(2)})\nno\\
&&=\res_{x_{0}}\frac{d}{dx}\left(z^{-1}
\delta\left(\frac{x^{-1}-x_{0}}{z}\right)\right)
\lambda(Y_{1}(e^{xL(1)}(-x^{-2})^{L(0)}v, x_{0})w_{(1)}\otimes w_{(2)})\nno\\
&&\quad +\res_{x_{0}}z^{-1}\delta\left(\frac{x^{-1}-x_{0}}{z}\right)
\frac{d}{dx}\lambda(Y_{1}(e^{xL(1)}(-x^{-2})^{L(0)}v, x_{0})w_{(1)}
\otimes w_{(2)})
\nonumber\\
&&\quad +\lambda(w_{(1)}\otimes \frac{d}{dx}Y_{2}^{*}(v, x)w_{(2)})\nno\\
&&=\res_{x_{0}}\frac{d}{dx}\left(z^{-1}
\delta\left(\frac{x^{-1}-x_{0}}{z}\right)\right)
\lambda(Y_{1}(e^{xL(1)}(-x^{-2})^{L(0)}v, x_{0})w_{(1)}\otimes w_{(2)})\nno\\
&&\quad +\res_{x_{0}}z^{-1}\delta\left(\frac{x^{-1}-x_{0}}{z}\right)
\lambda(Y_{1}(e^{xL(1)}L(1)(-x^{-2})^{L(0)}v, x_{0})w_{(1)}\otimes w_{(2)})
\nonumber\\
&&\quad -2\res_{x_{0}}z^{-1}\delta\left(\frac{x^{-1}-x_{0}}{z}\right)
\cdot\nno\\
&&\quad \quad\quad\quad \cdot
\lambda(Y_{1}(e^{xL(1)}L(0)x^{-1}(-x^{-2})^{L(0)}v, x_{0})w_{(1)}
\otimes w_{(2)})
\nonumber\\
&&\quad +\lambda(w_{(1)}\otimes Y_{2}^{*}(L(-1)v, x)w_{(2)}).
\end{eqnarray}
The first term
on the right-hand side of (13.9) is equal to
\begin{eqnarray}
\lefteqn{-\res_{x_{0}}x^{-2}\frac{d}{dx^{-1}}\left(z^{-1}
\delta\left(\frac{x^{-1}-x_{0}}{z}\right)\right)\cdot}\nno\\
&&\quad \quad\quad\quad \cdot
\lambda(Y_{1}(e^{xL(1)}(-x^{-2})^{L(0)}v, x_{0})w_{(1)}\otimes w_{(2)})\nno\\
&&=\res_{x_{0}}x^{-2}\frac{d}{dx_{0}}\left(z^{-1}
\delta\left(\frac{x^{-1}-x_{0}}{z}\right)\right)\cdot\nno\\
&&\quad \quad\quad\quad \cdot
\lambda(Y_{1}(e^{xL(1)}(-x^{-2})^{L(0)}v, x_{0})w_{(1)}\otimes w_{(2)})\nno\\
&&=-\res_{x_{0}}x^{-2}z^{-1}
\delta\left(\frac{x^{-1}-x_{0}}{z}\right)\cdot\nno\\
&&\quad \quad\quad\quad \cdot
\frac{d}{dx_{0}}\lambda(Y_{1}(e^{xL(1)}(-x^{-2})^{L(0)}v, x_{0})w_{(1)}
\otimes w_{(2)})\nno\\
&&=-\res_{x_{0}}x^{-2}z^{-1}
\delta\left(\frac{x^{-1}-x_{0}}{z}\right)\cdot\nno\\
&&\quad \quad\quad\quad \cdot
\lambda(Y_{1}(L(-1)e^{xL(1)}(-x^{-2})^{L(0)}v, x_{0})w_{(1)}
\otimes w_{(2)}),
\end{eqnarray}
where we have used the ``integration by parts'' property of $\res_{x_{0}}$.
By (13.10), (5.2.14) of \cite{FHL} and an appropriate analogue of (5.2.12)
of \cite{FHL},
the right-hand side
of (13.9) is equal to
\begin{eqnarray*}
&{\displaystyle \res_{x_{0}}z^{-1}\delta\left(\frac{x^{-1}-x_{0}}{z}\right)
\lambda(Y_{1}(e^{xL(1)}(-x^{-2})^{L(0)}L(-1)v, x_{0})w_{(1)}\otimes w_{(2)})}&
\nonumber\\
&{\displaystyle +\lambda(w_{(1)}\otimes Y_{2}^{*}(L(-1)v, x)w_{(2)})}&\nno\\
&{\displaystyle =(Y'_{P(z)}(L(-1)v, x)\lambda)(w_{(1)}\otimes w_{(2)}),}&
\end{eqnarray*}
proving the $L(-1)$-derivative property.
\epfv

Write
\begin{equation}
Y'_{P(z)}(\omega, x)=\sum_{n\in {\Bbb Z}}L'_{P(z)}(n)x^{-n-2}
\end{equation}
(recall that $\omega$ is the generator of the Virasoro algebra,
for our vertex operator algebra $(V, Y, {\bf 1}, \omega)$).
We call the eigenspaces of the operator $L'_{P(z)}(0)$ the {\it $P(z)$-weight
subspaces} or {\it
$P(z)$-homogeneous subspaces} of $(W_{1}\otimes W_{2})^{*}$, and we have the
corresponding notions of {\it $P(z)$-weight vector} (or {\it
$P(z)$-homogeneous
vector}) and  {\it $P(z)$-weight}.

We shall not discuss commutators $[Y'_{P(z)}(v_{1}, x_{1}),
Y'_{P(z)}(v_{2}, x_{2})]$ on $(W_{1}\otimes W_{2})^{*}$ (cf.
Proposition 5.2 of Part I), but we shall instead directly discuss the
``compatibility condition,'' which will lead to the Jacobi identity on
a suitable subspace of $(W_{1}\otimes W_{2})^{*}$.

Suppose that $G\in \mbox{\rm Hom}(W'_{3}, (W_{1}\otimes
W_{2})^{*})$ intertwines the two actions as in Proposition 13.1.
Then for $w'_{(3)}\in W'_{3}$,
$G(w'_{(3)})$ satisfies the following nontrivial and subtle condition
on $\lambda
\in (W_{1}\otimes W_{2})^{*}$: The formal Laurent series $Y'_{P(z)}(v,
x_{0})\lambda$ involves only finitely many negative powers of $x_{0}$
and
\begin{eqnarray}
\lefteqn{\tau_{P(z)}\left(x_{0}^{-1}\delta\left(\frac{x^{-1}_{1}-z}{x_{0}}
\right)
Y_{t}(v, x_{1})\right)\lambda=}\nno\\
&&=x_{0}^{-1}\delta\left(\frac{x^{-1}_{1}-z}{x_{0}}\right)
Y'_{P(z)}(v, x_{1})\lambda  \;\;\;\;\; \mbox{\rm for all}\;\;v\in V.
\end{eqnarray}
(Note that the two sides are not {\it a priori} equal for general
$\lambda\in (W_{1}\otimes W_{2})^{*}$.)
We call this the {\it $P(z)$-compatibility  condition} on
$\lambda\in (W_{1}\otimes W_{2})^{*}$.
(Note that this compatibility condition is different {}from the compatibility
condition for $\tau_{Q(z)}$.)

Let $W$ be a subspace of $(W_{1}\otimes W_{2})^{*}$.  We say that $W$
is {\it $P(z)$-compatible} if every element of $W$
satisfies the $P(z)$-compatibility condition. Also, we say that $W$ is
(${\Bbb C}$-){\it graded} (by $L'_{P(z)}(0)$)
if it is ${\Bbb C}$-graded by its $P(z)$-weight subspaces,
and that $W$ is a $V$-{\it module} (respectively, {\it
generalized $V$-module}) if $W$ is graded and is a module (respectively,
generalized module) when equipped with this grading and with the action of
$Y'_{P(z)}(\cdot, x)$. A sum of compatible
modules or generalized modules is clearly a generalized module. The
weight subspace of a subspace $W$ with weight $n\in {\Bbb C}$ will be
denoted $W_{(n)}$.

Given $G$ as above, it is clear that $G(W'_{3})$ is a $V$-module since $G$
intertwines the two actions of $V\otimes {\Bbb C}[t, t^{-1}]$. We have in
fact established that $G(W'_{3})$ is in addition a $P(z)$-compatible
$V$-module since
$G$ intertwines the full actions. Moreover, if
$H\in \mbox{\rm Hom}(W'_{3}, (W_{1}\otimes W_{2})^{*})$ intertwines
the two
actions of $V\otimes {\Bbb C}[t, t^{-1}]$, then $H$ intertwines the two actions
of $V\otimes \iota_{+}{\Bbb C}[t, t^{-1}, (z^{-1}-t)^{-1}]$ if and only if the
$V$-module $H(W'_{3})$ is $P(z)$-compatible.

Define
\begin{equation}
W_{1}\hboxtr_{P(z)}W_{2}=\sum_{W\in {\cal W}_{P(z)}}W =\bigcup_{W\in
{\cal W}_{P(z)}} W\subset
(W_{1}\otimes W_{2})^{*},
\end{equation}
where ${\cal W}_{P(z)}$ is the set of all $P(z)$-compatible
modules in $(W_{1}\otimes W_{2})^{*}$.
Then $W_{1}\hboxtr_{P(z)}W_{2}$
is a $P(z)$-compatible generalized module and coincides with the sum (or
union) of the
images $G(W'_{3})$ of modules $W'_{3}$ under the maps $G$ as above.
Moreover, for any $V$-module $W_{3}$ and any map
$H: W'_{3}\to W_{1}\hboxtr_{P(z)} W_{2}$ of generalized modules,
$H(W'_{3})$ is $P(z)$-compatible and hence $H$ intertwines the two actions
of $V\otimes \iota_{+}{\Bbb C}[t, t^{-1}, (z^{-1}-t)^{-1}]$. Thus we have:
\begin{propo}
The subspace $W_{1}\hboxtr_{P(z)}W_{2}$ of $(W_{1}\otimes W_{2})^{*}$ is
a generalized module with the following
property: Given  any $V$-module $W_{3}$,
there is a natural linear isomorphism determined by (13.3) between the space of
all $P(z)$-intertwining maps of type ${W_{3}}\choose {W_{1}\;W_{2}}$
and the space of all maps of generalized modules {from} $W'_{3}$ to
$W_{1}\hboxtr _{P(z)}W_{2}$.\epf
\end{propo}

For rational vertex operator algebras, we have the following straightforward
result, proved exactly as in Proposition 5.6 of Part I:

\begin{propo}
Let $V$ be a rational vertex operator algebra and $W_{1}$, $W_{2}$ two
$V$-modules. Then the generalized module
$W_{1}\hboxtr _{P(z)}W_{2}$ is a module.\epf
\end{propo}

Now we assume that $W_{1}\hboxtr_{P(z)} W_{2}$ is a module (which
occurs if $V$ is rational, by
the last proposition).  In this case, we define a $V$-module
$W_{1}\boxtimes_{P(z)} W_{2}$ by
\begin{equation}
W_{1}\boxtimes_{P(z)} W_{2}=(W_{1}\hboxtr_{P(z)}W_{2})'
\end{equation}
 and we write the corresponding
action as $Y_{P(z)}$. Applying Proposition 13.4 to the special module
$W_{3}=W_{1}\boxtimes_{P(z)} W_{2}$ and the identity map $W'_{3}\to
W_{1}\hboxtr_{P(z)} W_{2}$, we obtain using
(13.3) a canonical $P(z)$-intertwining map of type
${W_{1}\boxtimes_{P(z)} W_{2}}\choose {W_{1}W_{2}}$, which we denote
\begin{eqnarray}
\boxtimes_{P(z)}: W_{1}\otimes W_{2}&\to &
\overline{W_{1}\boxtimes_{P(z)} W_{2}}\nno\\
w_{(1)}\otimes  w_{(2)}&\mapsto& w_{(1)} \boxtimes_{P(z)}w_{(2)}.
\end{eqnarray}

It is easy to verify, just as in Proposition 5.7 of Part I:

\begin{propo}
The $P(z)$-product $(W_{1}\boxtimes_{P(z)} W_{2}, Y_{P(z)};
\boxtimes_{P(z)})$
is a $P(z)$-tensor product of $W_{1}$ and $W_{2}$.\epf
\end{propo}

More generally, dropping the assumption that $W_{1}\hboxtr_{P(z)}
W_{2}$ is a module, we have, by imitating the proof of Proposition 5.8
of Part I:

\begin{propo}
The $P(z)$-tensor product of $W_{1}$ and $W_{2}$ exists (and is given
by (13.14) if and only if $W_{1}\hboxtr_{P(z)} W_{2}$ is a module.\epf
\end{propo}

We observe that any element of
$W_{1}\hboxtr_{P(z)} W_{2}$ is an element $\lambda$
of $(W_{1}\otimes W_{2})^{*}$ satisfying:

\begin{description}
\item[The $P(z)$-compatibility condition
{\rm (recall (13.12))}] \hfill

{\bf (a)} The  {\it $P(z)$-lower
truncation condition}:
For all $v\in V$, the formal Laurent series $Y'_{P(z)}(v, x)\lambda$
involves only finitely many negative
powers of $x$.

{\bf (b)} The following formula holds:
\begin{eqnarray}
\lefteqn{\tau_{P(z)}\left(x_{0}^{-1}\delta\left(\frac{x^{-1}_{1}-z}{x_{0}}
\right)
Y_{t}(v, x_{1})\right)\lambda=}\nno\\
&&=x_{0}^{-1}\delta\left(\frac{x^{-1}_{1}-z}{x_{0}}\right)
Y'_{P(z)}(v, x_{1})\lambda  \;\;\;\;\; \mbox{\rm for all}\;\;v\in V.
\end{eqnarray}

\item[The $P(z)$-local grading-restriction  condition]\hfill

{\bf (a)} The {\it $P(z)$-grading condition}:
$\lambda$ is a (finite) sum of
weight vectors of $(W_{1}\otimes W_{2})^{*}$.

{\bf (b)} Let $W_{\lambda}$ be the smallest subspace of $(W_{1}\otimes
W_{2})^{*}$ containing $\lambda$ and stable under the component
operators $\tau_{P(z)}(v\otimes t^{n})$ of the operators $Y'_{P(z)}(v,
x)$ for $v\in V$, $n\in {\Bbb Z}$. Then the weight spaces
$(W_{\lambda})_{(n)}$, $n\in {\Bbb C}$, of the (graded) space
$W_{\lambda}$ have the properties
\begin{eqnarray}
&\mbox{\rm dim}\ (W_{\lambda})_{(n)}<\infty \;\;\;\mbox{\rm for}\
n\in {\Bbb C},&\\
&(W_{\lambda})_{(n)}=0 \;\;\;\mbox{\rm for $n$ whose real part is
sufficiently small.}&
\end{eqnarray}
\end{description}

We shall call the compatibility condition, the lower truncation condition,
the local grading-restriction condition and the grading condition
for $\tau_{Q(z)}$ the
{\it $Q(z)$-compatibility condition},
the {\it $Q(z)$-lower truncation condition},
the {\it $Q(z)$-local grading-restriction condition} and the {\it
$Q(z)$-grading condition} when it is necessary to distinguish
those conditions for
$\tau_{P(z)}$ and for $\tau_{Q(z)}$.

The next lemma allows us to
use the results proved in the construction of the $Q(z^{-1})$-tensor
product to give another, much more useful,
construction of the $P(z)$-tensor product.
Let $\psi: W_{1}\otimes W_{2}\to  W_{1}\otimes W_{2}$ be the linear map
defined by
\begin{equation}
\psi(w_{(1)}\otimes w_{(2)})=e^{-z^{-1}L(1)}w_{(1)}\otimes
e^{(-2\log z+\pi i)L(0)}e^{z^{-1}L(1)}w_{(2)}
\end{equation}
where $w_{(1)}\in W_{1}$ and $w_{(2)}\in W_{2}$ and let $\psi^{*}:
(W_{1}\otimes W_{2})^{*}\to (W_{1}\otimes W_{2})^{*}$ be the adjoint
of $\psi$.  Note that $\psi$ is a linear isomorphism and hence so is
$\psi^{*}$.

\begin{lemma}
For any $f\in (W_{1}\otimes W_{2})^{*}$, we have
\begin{eqnarray}
\lefteqn{\tau_{P(z)}\left(x_{0}^{-1}
\delta\left(\frac{x_{1}^{-1}-z}{x_{0}}\right)
Y_{t}(v, x_{1})\right)\psi^{*}(f)=}\nno\\
&&=(zx_{0})^{-1}\psi^{*}\biggl(\tau_{Q(z^{-1})}\biggl(zx_{0}x_{1}\delta
\left(\frac{z^{-1}+x_{0}^{-1}}{(zx_{0}x_{1})^{-1}}\right)\cdot\nno\\
&&\hspace{6em}\cdot Y_{t}(e^{zx_{0}x_{1}L(1)}(x_{0}x_{1})^{-2L(0)}v,
x_{0}^{-1})\biggr)f\biggr).
\end{eqnarray}
\end{lemma}
\pf
For any $w_{(1)}\in W_{1}$ and $w_{(2)}\in W_{2}$, by the definitions
of $\tau_{P(z)}$ and $\psi$ and by the properties of the formal
$\delta$-functions, we have
\begin{eqnarray}
\lefteqn{\left(\tau_{P(z)}\left(x_{0}^{-1}\delta\left(\frac{x_{1}^{-1}-z}
{x_{0}}\right)Y_{t}(v, x_{1})\right)\psi^{*}(f)\right)(w_{(1)}
\otimes w_{(2)})=}\nno\\
&&=z^{-1}\delta\left(\frac{x^{-1}_{1}-x_{0}}{z}\right)
(\psi^{*}(f))(Y_{1}(e^{x_{1}L(1)}(-x_{1}^{2})^{-L(0)}v,
x_{0})w_{(1)}\otimes w_{(2)})\nonumber\\
&&\hspace{1em}+x^{-1}_{0}\delta\left(\frac{z-x^{-1}_{1}}{-x_{0}}\right)
(\psi^{*}(f))(w_{(1)}\otimes Y_{2}^{*}(v, x_{1})w_{(2)})\nno\\
&&=z^{-1}\delta\left(\frac{x^{-1}_{1}-x_{0}}{z}\right)
f(e^{-z^{-1}L(1)}Y_{1}(e^{x_{1}L(1)}(-x_{1}^{2})^{-L(0)}v,
x_{0})w_{(1)}\otimes\nno\\
&&\hspace{6em}\otimes e^{(-2\log z+\pi i)L(0)}e^{z^{-1}L(1)}w_{(2)})\nonumber\\
&&\hspace{1em}+x^{-1}_{0}\delta\left(\frac{z-x^{-1}_{1}}{-x_{0}}\right)
f(e^{-z^{-1}L(1)}w_{(1)}\otimes\nno\\
&&\hspace{6em}\otimes  e^{(-2\log z+\pi i)L(0)}e^{z^{-1}L(1)}
Y_{2}^{*}(v, x_{1})w_{(2)}).
\end{eqnarray}
As in the proof of Lemma 5.2.3 of \cite{FHL}, together with the use of
the adjoints, acting on $W'_{2}$, of the operators $Y^{*}_{2}(v, x)$, $L(0)$
and $L(1)$ acting on $W_{2}$ (recall (3.22) in Part I and (5.2.10) in
\cite{FHL}), we have
\begin{eqnarray}
e^{\zeta L(1)}Y_{1}(v, x)e^{-\zeta L(1)}&=&Y_{1}(
e^{\zeta(1-\zeta x)L(1)}(1-\zeta x)^{-2L(0)}v,
x/(1-\zeta x)),\;\;\;\;\;\;\;\;\;\;\\
e^{\zeta L(1)}Y^{*}_{2}(v, x)e^{-\zeta L(1)}&=&Y_{2}^{*}(v, x-\zeta),\\
e^{\zeta L(0)}Y_{2}^{*}(v, x)e^{-\zeta L(0)}
&=&Y^{*}(e^{-\zeta L(0)}v, e^{-\zeta} x)
\end{eqnarray}
for any complex number $\zeta$.
Using (13.22)--(13.24), we see that the right-hand side of (13.21) becomes
\begin{eqnarray}
\lefteqn{\dps z^{-1}\delta\left(\frac{x^{-1}_{1}-x_{0}}{z}\right)
f(Y_{1}(e^{-z^{-1}(1+z^{-1}x_{0})L(1)}(1+z^{-1}x_{0})^{-2L(0)}\cdot} \nno\\
&&\hspace{2em}\cdot e^{x_{1}L(1)}(-x_{1}^{2})^{-L(0)}v,
x_{0}/(1+z^{-1}x_{0}))e^{-z^{-1}L(1)}w_{(1)}\otimes\nno\\
&&\hspace{4em}\otimes e^{(-2\log z+\pi i)L(0)}e^{z^{-1}L(1)}w_{(2)})\nonumber\\
&& +x^{-1}_{0}\delta\left(\frac{z-x^{-1}_{1}}{-x_{0}}\right)
f(e^{-z^{-1}L(1)}w_{(1)}\otimes \nno\\
&&\hspace{2em}\otimes Y_{2}^{*}((-z^{2})^{L(0)}v,
-z^{2}(x_{1}-z^{-1}))
e^{(-2\log z+\pi i)L(0)}e^{z^{-1}L(1)}w_{(2)}).\hspace{3em}
\end{eqnarray}
Using the properties of the formal $\delta$-functions, (13.25) is equal to
\begin{eqnarray}
\lefteqn{ x_{1}\delta\left(\frac{z+x_{0}}{x_{1}^{-1}}\right)
f(Y_{1}(e^{-z^{-2}x_{1}^{-1}L(1)}(zx_{1})^{2L(0)}
e^{x_{1}L(1)}\cdot} \nno\\
&&\hspace{2em}\cdot (-x_{1}^{2})^{-L(0)}v,
zx_{0}x_{1})e^{-z^{-1}L(1)}w_{(1)}\otimes
e^{(-2\log z+\pi i)L(0)}e^{z^{-1}L(1)}w_{(2)})\nonumber\\
&&+x^{-1}_{0}\delta\left(\frac{z-x^{-1}_{1}}{-x_{0}}\right)
f(e^{-z^{-1}L(1)}w_{(1)}\otimes \nno\\
&&\hspace{2em}\otimes
Y_{2}^{*}((-z^{2})^{L(0)}v,
zx_{0}x_{1})e^{(-2\log z+\pi i)L(0)}e^{z^{-1}L(1)}w_{(2)}).
\end{eqnarray}
By the formula
\begin{equation}
x^{L(0)}e^{x_{1}L(1)}x^{-L(0)}=e^{(x_{1}/x)L(1)}
\end{equation}
(cf. (5.3.1)--(5.3.3) in \cite{FHL}) and the definitions
of $Y_{1}^{*}$ and $Y_{2}^{*}$ and of $\tau_{Q(z^{-1})}$, (13.26) is equal to
\begin{eqnarray}
\lefteqn{ x_{1}\delta\left(\frac{z+x_{0}}{x_{1}^{-1}}\right)
f(Y_{1}((-z^{2})^{L(0)}v,
zx_{0}x_{1})e^{-z^{-1}L(1)}w_{(1)}\otimes }\nno\\
&&\hspace{2em}\otimes e^{(-2\log z+\pi i)L(0)}e^{z^{-1}L(1)}w_{(2)})\nonumber\\
&&\hspace{1em}+x^{-1}_{0}\delta\left(\frac{z-x^{-1}_{1}}{-x_{0}}\right)
f(e^{-z^{-1}L(1)}w_{(1)}\otimes Y_{2}(e^{zx_{0}x_{1}L(1)}\cdot \nno\\
&&\hspace{2em}\cdot (-(zx_{0}x_{1})^{-2})^{L(0)}(-z^{2})^{L(0)}v,
(zx_{0}x_{1})^{-1})
e^{(-2\log z+\pi i)L(0)}e^{z^{-1}L(1)}w_{(2)})\nno\\
&& =z^{-1}\delta\left(\frac{x^{-1}_{1}-x_{0}}{z}\right)
f(Y_{1}^{*}((-(zx_{0}x_{1})^{-2})^{L(0)}e^{-(zx_{0}x_{1})^{-1}L(1)}\cdot \nno\\
&&\hspace{2em}\cdot (-z^{2})^{L(0)}v,
(zx_{0}x_{1})^{-1})e^{-z^{-1}L(1)}w_{(1)}\otimes
e^{(-2\log z+\pi i)L(0)}e^{z^{-1}L(1)}w_{(2)})\nonumber\\
&&\hspace{1em} -z^{-1}\delta\left(\frac{x_{0}-
x_{1}^{-1}}{-z}\right)
f(e^{-z^{-1}L(1)}w_{(1)}\otimes Y_{2}(e^{zx_{0}x_{1}L(1)}\cdot \nno\\
&&\hspace{2em}\cdot (x_{0}x_{1})^{-2L(0)}v,
(zx_{0}x_{1})^{-1})
e^{(-2\log z+\pi i)L(0)}e^{z^{-1}L(1)}w_{(2)})\nno\\
&& =(zx_{0})^{-1}x_{0}\delta\left(\frac{(zx_{0}x_{1})^{-1}-z^{-1}}
{x_{0}^{-1}}\right)
f(Y_{1}^{*}(e^{zx_{0}x_{1}L(1)}
(x_{0}x_{1})^{-2L(0)}v, \nno\\
&&\hspace{2em}(zx_{0}x_{1})^{-1}) e^{-z^{-1}L(1)}w_{(1)}\otimes
e^{(-2\log z+\pi i)L(0)}e^{z^{-1}L(1)}w_{(2)})\nonumber\\
&&\hspace{1em} -(zx_{0})^{-1}x_{0}\delta\left(\frac{z^{-1}-
(zx_{0}x_{1})^{-1}}{-x_{0}^{-1}}\right)
f(e^{-z^{-1}L(1)}w_{(1)}\otimes \nno\\
&&\hspace{2em}\otimes Y_{2}(e^{zx_{0}x_{1}L(1)}
(x_{0}x_{1})^{-2L(0)}v,
(zx_{0}x_{1})^{-1})
e^{(-2\log z+\pi i)L(0)}e^{z^{-1}L(1)}w_{(2)})\nno\\
&& =(zx_{0})^{-1}\psi^{*}\biggl(\tau_{Q(z^{-1})}\biggl(zx_{0}x_{1}\delta
\left(\frac{z^{-1}+x_{0}^{-1}}{(zx_{0}x_{1})^{-1}}\right)\cdot \nno\\
&&\hspace{2em}\cdot Y_{t}(e^{zx_{0}x_{1}L(1)}
(x_{0}x_{1})^{-2L(0)}v, x_{0}^{-1})\biggr)f\biggr)
(w_{(1)}\otimes w_{(2)}),
\end{eqnarray}
proving (13.20).
\epfv

The following result includes the analogues of Theorem 6.1 and
Proposition 6.2 of \cite{HL1}, which are parts of Theorem 6.3 of
\cite {HL1} (stated as Theorem I.15 above).  We state the result
as one theorem since the proofs are intermeshed with one another.
\begin{theo}
An element $f\in (W_{1}\otimes W_{2})^{*}$ satisfies the
$P(z)$-compatibility condition if and only if $(\psi^{*})^{-1}(f)$
satisfies the $Q(z^{-1})$-compatibility condition. In this case, $f$
satisfies the $P(z)$-local grading-restriction condition if and only
if $(\psi^{*})^{-1}(f)$ satisfies the $Q(z^{-1})$-local
grading-restriction condition. If $f$ satisfies both conditions, then
$\tau_{P(z)}(v\otimes t^{n})f$, $v\in V$, $n\in {\Bbb Z}$, also
satisfies both conditions and we have the Jacobi identity acting on
$f$:
\begin{eqnarray}
\lefteqn{x_{0}^{-1}\delta
\left({\displaystyle\frac{x_{1}-x_{2}}{x_{0}}}\right)
Y'_{P(z)}(u, x_{1})
Y'_{P(z)}(v, x_{2})f}\nno\\
&&\hspace{2ex}-x_{0}^{-1} \delta
\left({\displaystyle\frac{x_{2}-x_{1}}{-x_{0}}}\right)
Y'_{P(z)}(v, x_{2})
Y'_{P(z)}(u, x_{1})f\nonumber \\
&&=x_{2}^{-1} \delta
\left({\displaystyle\frac{x_{1}-x_{0}}{x_{2}}}\right)
Y'_{P(z)}(Y(u, x_{0})v,
x_{2})f.
\end{eqnarray}
\end{theo}
\pf
Let $f$ be an element of $(W_{1}\otimes W_{2})^{*}$ satisfying
the $P(z)$-compatibility
condition.  Since $f$ satisfies
the $P(z)$-lower truncation condition,
$$Y'_{P(z)}((1+zx_{0}^{-1})^{-2L(0)}v,
x_{0}^{-1}(1+zx_{0}^{-1})^{-1})f$$
is well defined for any $v\in V$ and is a Laurent series in $x_{0}^{-1}$
containing
 only finitely
many negative powers of $x_{0}^{-1}$.
{}From the
$P(z)$-compatibility condition  and the
fundamental property of the $\delta$-function, we have
\begin{eqnarray}
\lefteqn{\tau_{P(z)}\left(x_{0}^{-1}
\delta\left(\frac{x_{1}^{-1}-z}{x_{0}}\right)
Y_{t}((x_{0}x_{1})^{2L(0)}v,
x_{1})\right)f}\nno\\
&&=x_{0}^{-1}\delta\left(\frac{x^{-1}_{1}-z}{x_{0}}\right)
Y'_{P(z)}((x_{0}x_{1})^{2L(0)}v, x_{1})f\nno\\
&&=x_{0}^{-1}\delta\left(\frac{x^{-1}_{1}-z}{x_{0}}\right)
Y'_{P(z)}((1+zx_{0}^{-1})^{-2L(0)}v,
x_{0}^{-1}(1+zx_{0}^{-1})^{-1})f.\nno\\
&&
\end{eqnarray}
Substituting $e^{-zx_{0}x_{1}L(1)}v$ for $v$ and again using the
fundamental property of the $\delta$-function on the right-hand side,
we find that
\begin{eqnarray}
\lefteqn{\tau_{P(z)}\left(x_{0}^{-1}
\delta\left(\frac{x_{1}^{-1}-z}{x_{0}}\right)
Y_{t}((x_{0}x_{1})^{2L(0)}e^{-zx_{0}x_{1}L(1)}v,
x_{1})\right)f=}\nno\\
&&=x_{0}^{-1}\delta\left(\frac{x^{-1}_{1}-z}{x_{0}}\right)
Y'_{P(z)}((1+zx_{0}^{-1})^{-2L(0)}\cdot \nno\\
&&\hspace{6em}\cdot e^{-z(1+zx_{0}^{-1})^{-1}L(1)}v,
x_{0}^{-1}(1+zx_{0}^{-1})^{-1})f.
\end{eqnarray}

{}From Lemma 13.8, we have
\begin{eqnarray}
\lefteqn{\tau_{P(z)}\left(x_{0}^{-1}
\delta\left(\frac{x_{1}^{-1}-z}{x_{0}}\right)
Y_{t}((x_{0}x_{1})^{2L(0)}e^{-zx_{0}x_{1}L(1)}v,
x_{1})\right)f=}\nno\\
&&=(zx_{0})^{-1}\psi^{*}\left(\tau_{Q(z^{-1})}\left(zx_{0}x_{1}\delta
\left(\frac{z^{-1}+x_{0}^{-1}}{(zx_{0}x_{1})^{-1}}\right)
Y_{t}(v, x_{0}^{-1})\right)
(\psi^{*})^{-1}(f)\right)\nno\\
&&
\end{eqnarray}
and we obtain
\begin{eqnarray}
\lefteqn{x_{0}^{-1}\delta\left(\frac{x^{-1}_{1}-z}{x_{0}}\right)
Y'_{P(z)}((1+zx_{0}^{-1})^{-2L(0)}\cdot} \nno\\
&&\hspace{6em}\cdot e^{-z(1+zx_{0}^{-1})^{-1}L(1)}v,
x_{0}^{-1}(1+zx_{0}^{-1})^{-1})f\nno\\
&&=(zx_{0})^{-1}\psi^{*}\left(\tau_{Q(z^{-1})}\left(zx_{0}x_{1}\delta
\left(\frac{z^{-1}+x_{0}^{-1}}{(zx_{0}x_{1})^{-1}}\right)
Y_{t}(v, x_{0}^{-1})\right)
(\psi^{*})^{-1}(f)\right).\nno\\
&&
\end{eqnarray}
Extracting $\res_{x_{1}^{-1}}$ gives
\begin{eqnarray}
\lefteqn{Y'_{P(z)}((1+zx_{0}^{-1})^{-2L(0)}e^{-z(1+zx_{0}^{-1})^{-1}L(1)}v,
x_{0}^{-1}(1+zx_{0}^{-1})^{-1})f}\nno\\
&&=\res_{x_{1}^{-1}}(zx_{0})^{-1}\cdot\nno\\
&&\hspace{4em}\cdot \psi^{*}\left(\tau_{Q(z^{-1})}
\left(zx_{0}x_{1}\delta
\left(\frac{z^{-1}+x_{0}^{-1}}{(zx_{0}x_{1})^{-1}}\right)
Y_{t}(v, x_{0}^{-1})\right)
(\psi^{*})^{-1}(f)\right)\nno\\
&&=\res_{y}\psi^{*}\left(\tau_{Q(z^{-1})}
\left(y^{-1}\delta
\left(\frac{z^{-1}+x_{0}^{-1}}{y}\right)
Y_{t}(v, x_{0}^{-1})\right)
(\psi^{*})^{-1}(f)\right)\nno\\
&&=\psi^{*}(Y'_{Q(z^{-1})}(v, x^{-1}_{0})(\psi^{*})^{-1}(f)).
\end{eqnarray}
Since the left-hand side of (13.34) is a well-defined formal
Laurent series in $x_{0}^{-1}$ containing
only finitely many negative powers of $x_{0}^{-1}$, the right-hand side
of (13.34) is also such a Laurent series. Thus $(\psi^{*})^{-1}(f)$ satisfies
the
$Q(z^{-1})$-lower truncation condition. Using (13.33) and (13.34), we have
\begin{eqnarray}
\lefteqn{(zx_{0})^{-1}\psi^{*}\left(\tau_{Q(z^{-1})}\left(zx_{0}x_{1}\delta
\left(\frac{z^{-1}+x_{0}^{-1}}{(zx_{0}x_{1})^{-1}}\right)
Y_{t}(v, x_{0}^{-1})\right)
(\psi^{*})^{-1}(f)\right)}\nno\\
&&=x_{1}\delta
\left(\frac{x_{0}+z}{x_{1}^{-1}}\right)\psi^{*}\left(
Y'_{Q(z^{-1})}(v, x_{0}^{-1})
(\psi^{*})^{-1}(f)\right),\nno\\
&&
\end{eqnarray}
or equivalently, replacing $x_0$ by $x_{0}^{-1}$ and then replacing
$x_1$ by $z^{-1}x_{0}x_{1}^{-1}$ and applying $(\psi^{*})^{-1}$,
\begin{eqnarray}
\lefteqn{\tau_{Q(z^{-1})}\left(x_{1}^{-1}\delta
\left(\frac{z^{-1}+x_{0}}{x_{1}}\right)
Y_{t}(v, x_{0})\right)
(\psi^{*})^{-1}(f)=}\nno\\
&&\hspace{1em}=x_{1}^{-1}\delta
\left(\frac{z^{-1}+x_{0}}{x_{1}}\right)
Y'_{Q(z^{-1})}(v, x_{0})
(\psi^{*})^{-1}(f).
\end{eqnarray}
Thus $(\psi^{*})^{-1}(f)$ satisfies the
$Q(z^{-1})$-compatibility condition.

Conversely, if $(\psi^{*})^{-1}(f)$ satisfies the
$Q(z^{-1})$-compatibility condition, an analogue of the proof above
shows that $f$ satisfies the $P(z)$-compatibility condition: The
$Q(z^{-1})$-compatibility condition (13.36) and Lemma 13.8 give
analogues of (13.33) and (13.34), and these analogues imply that $f$
satisfies the $P(z)$-compatibility condition. This finishes the proof
of the first part of the theorem.

Now assume that $f$ satisfies the $P(z)$-compatibility condition or
equivalently, that $(\psi^{*})^{-1}(f)$ satisfies the
$Q(z^{-1})$-compatibility condition.  By Proposition 5.1 and Theorem
6.1 of \cite{HL1}, which are parts of
Theorem \ref{6-3} above (Theorem 6.3 of
\cite{HL1}), we have the $L(-1)$-derivative property
\begin{equation}
\frac{d}{dx}Y'_{Q(z^{-1})}(v, x)(\psi^{*})^{-1}(f)
=Y'_{Q(z^{-1})}(L(-1)v, x)(\psi^{*})^{-1}(f)
\end{equation}
(and this in fact holds for any $g\in (W_{1}\otimes W_{2})^{*}$ in
place of $(\psi^{*})^{-1}(f)$) and the Jacobi identity
\begin{eqnarray}
\lefteqn{x_{0}^{-1}\delta
\left({\displaystyle\frac{x_{1}-x_{2}}{x_{0}}}\right)Y'_{Q(z^{-1})}(u, x_{1})
Y'_{Q(z^{-1})}(v, x_{2})(\psi^{*})^{-1}(f)}\nno\\
&&\hspace{2ex}-x_{0}^{-1} \delta
\left({\displaystyle\frac{x_{2}-x_{1}}{-x_{0}}}\right)Y'_{Q(z^{-1})}(v, x_{2})
Y'_{Q(z^{-1})}(u, x_{1})(\psi^{*})^{-1}(f)\nonumber \\
&&=x_{2}^{-1} \delta
\left({\displaystyle\frac{x_{1}-x_{0}}{x_{2}}}\right)Y'_{Q(z^{-1})}(Y(u,
x_{0})v,
x_{2})(\psi^{*})^{-1}(f).
\end{eqnarray}
for $u,v \in V$.
Since the Virasoro element $\omega=L(-2){\bf 1}$ is quasi-primary,
that is, $L(1)\omega=0$, and its weight is $2$,
(13.34) with $v=\omega$ gives
\begin{eqnarray}\label{1-0}
\lefteqn{\sum_{n\in {\Bbb Z}}(L'_{P(z)}(n)f)x^{n+2}_{0}
(1+zx_{0}^{-1})^{n-2}=}\nno\\
&&=Y'_{P(z)}((1+zx_{0}^{-1})^{-2L(0)}e^{-z(1+zx_{0}^{-1})^{-1}L(1)}\omega,
x_{0}^{-1}(1+zx_{0}^{-1})^{-1})f\nno\\
&&=\psi^{*}(Y'_{Q(z^{-1})}(\omega, x_{0}^{-1})(\psi^{*})^{-1}(f))\nno\\
&&=\sum_{n\in {\Bbb Z}}\psi^{*}(L'_{Q(z^{-1})}(n)(\psi^{*})^{-1}(f))
x_{0}^{n+2}.
\end{eqnarray}
where $L'_{P(z)}(n)$, $n\in \Bbb{Z}$, are defined in (13.11) and
$L'_{Q(z^{-1})}(n)$, $n\in \Bbb{Z}$, are defined analogously as the
coefficients of the vertex operator $Y'_{Q(z^{-1})}(\omega, x)$.
Taking the coefficient of $x^{3}_{0}$ on both sides of (\ref{1-0})
and noting
that for $n\in \Bbb{Z}$, the coefficient of $x^{3}_{0}$ in $x^{n+2}_{0}
(1+zx_{0}^{-1})^{n-2}$ is $\delta_{n, 1}$, we obtain
$$L'_{P(z)}(1)f=\psi^{*}(L'_{Q(z^{-1})}(1)(\psi^{*})^{-1}(f))$$
or equivalently,
\begin{equation}
(\psi^{*})^{-1}(L'_{P(z)}(1)f)=L'_{Q(z^{-1})}(1)(\psi^{*})^{-1}(f).
\end{equation}
Similarly, taking the coefficient of $x^{2}_{0}$ on both sides of
(\ref{1-0}) and using (13.40), we have
\begin{equation}
(\psi^{*})^{-1}(L'_{P(z)}(0)f)=(L'_{Q(z^{-1})}(0)+zL'_{Q(z^{-1})}(1))
(\psi^{*})^{-1}(f).
\end{equation}
{}From the Jacobi identity (13.38) together with (13.37), we obtain
the usual commutator formulas for
$$[L'_{Q(z^{-1})}(m),
L'_{Q(z^{-1})}(n)]$$ and for
$$[L'_{Q(z^{-1})}(m), Y'_{Q(z^{-1})}(v,
x)]$$
($m, n\in {\Bbb Z}$, $v\in V$), acting on $(\psi^{*})^{-1}(f)$.
In particular, we have
\begin{equation}
[L'_{Q(z^{-1})}(0), L'_{Q(z^{-1})}(1)](\psi^{*})^{-1}(f)
=-L'_{Q(z^{-1})}(1)(\psi^{*})^{-1}(f)
\end{equation}
and
\begin{eqnarray}
\lefteqn{[L'_{Q(z^{-1})}(1), Y'_{Q(z^{-1})}(v, x)](\psi^{*})^{-1}(f)=}\nno\\
&&=Y'_{Q(z^{-1})}((L(1)+2xL(0)+x^{2}L(-1))v, x)(\psi^{*})^{-1}(f).
\end{eqnarray}

By Proposition 6.2 of \cite{HL1} (part of Theorem \ref{6-3} above),
the subspace of $(W_{1}\otimes W_{2})^{*}$ consisting of the elements
satisfying the $Q(z^{-1})$-compatibility condition is stable under the
operators $\tau_{Q(z^{-1})}(v\otimes t^{n})$ for $v\in V$ and $n\in {\Bbb
Z}$, and in particular, under the operators $L'_{Q(z^{-1})}(n)$.
Thus by (13.40), (13.41) and what has been proved above, the subspace of
$(W_{1}\otimes W_{2})^{*}$ consisting of the elements satisfying the
$P(z)$-compatibility condition is stable under the operator
$L'_{P(z)}(1)$ and $L'_{P(z)}(0)$. Also, by (13.40)--(13.42), we have
\begin{equation}
[L'_{P(z)}(0), L'_{P(z)}(1)]f=-L'_{P(z)}(1)f.
\end{equation}

We also see that $\tau_{P(z)}(v\otimes t^{n})$ preserves the space of
elements of $(W_{1}\otimes W_{2})^{*}$ satisfying the
$P(z)$-compatibility condition. Indeed, replacing $x_{0}^{-1}$ in
(13.34) by $x(1-zx)^{-1}$, we have
\begin{eqnarray}
\lefteqn{Y'_{P(z)}((1-zx)^{2L(0)}e^{-z(1-zx)L(1)}v,
x)f=}\nno\\
&&=\psi^{*}(Y'_{Q(z^{-1})}(v, x(1-zx)^{-1})(\psi^{*})^{-1}(f)),
\end{eqnarray}
or equivalently,
\begin{eqnarray}
\lefteqn{Y'_{P(z)}(v,
x)f=}\nno\\
&&=\psi^{*}(Y'_{Q(z^{-1})}(e^{z(1-zx)L(1)}(1-zx)^{-2L(0)}v,
x(1-zx)^{-1})(\psi^{*})^{-1}(f)),\nno\\
&&
\end{eqnarray}
and we now invoke the already-established
equivalence between the compatibility conditions for $f$ and
$(\psi^{*})^{-1}(f)$.

Now suppose that $f$ satisfies the $P(z)$-compatibility condition and
that $L'_{P(z)}(1)$ acts nilpotently on $f$ (that is,
$(L'_{P(z)}(1))^{n}f=0$ for large $n$). Then of course $L'_{P(z)}(1)$
also acts nilpotently on such elements as $e^{-zL'_{P(z)}(1)}f$. By
(13.40) and the comments above,
$e^{-zL'_{Q(z^{-1})}(1)}(\psi^{*})^{-1}(f)$ is also well defined (in
the analogous sense).
{}From (13.37) and (13.43), we obtain, as in the proof of (5.2.38) of
\cite{FHL} (which is valid here since we are using the two formal
variables $x_{0}$ and $x$ and since (13.43) remains valid with
$(\psi^{*})^{-1}(f)$ replaced by $L'_{Q(z^{-1})}(1)(\psi^{*})^{-1}(f)$),
\begin{eqnarray}
\lefteqn{e^{x_{0}L'_{Q(z^{-1})}(1)}Y'_{Q(z^{-1})}(v, x)
e^{-x_{0}L'_{Q(z^{-1})}(1)}
(\psi^{*})^{-1}(f)=}\nno\\
&&=Y'_{Q(z^{-1})}(e^{x_{0}(1-x_{0}x)L(1)}(1-x_{0}x)^{-2L(0)}v,
x(1-x_{0}x)^{-1})
(\psi^{*})^{-1}(f).\nno\\
&&
\end{eqnarray}

The coefficient of each monomial in $x$ on the right-hand side of
(13.47) is a (terminating) polynomial in $x_{0}$, and so the same is
true of the left-hand side. Thus we may substitute $z$ for $x_{0}$
and we obtain
\begin{eqnarray}
\lefteqn{e^{zL'_{Q(z^{-1})}(1)}Y'_{Q(z^{-1})}(v, x)
e^{-zL'_{Q(z^{-1})}(1)}(\psi^{*})^{-1}(f)=}
\nno\\
&&=Y'_{Q(z^{-1})}(e^{z(1-zx)L(1)}(1-zx)^{-2L(0)}v, x(1-zx)^{-1})
(\psi^{*})^{-1}(f).\nno\\
&&
\end{eqnarray}

Combining (13.46) and (13.48), we have the following formula, which
relates $Y'_{P(z)}(v, x)$ and $Y'_{Q(z^{-1})}(v, x)$ by conjugation,
when acting on $f$:
\begin{equation}
Y'_{P(z)}(v,
x)f=\psi^{*}(e^{zL'_{Q(z^{-1})}(1)}Y'_{Q(z^{-1})}(v, x)
e^{-zL'_{Q(z^{-1})}(1)}
(\psi^{*})^{-1}(f)).
\end{equation}
(Recall that we have just seen that $L'_{Q(z^{-1})}(1)$ acts nilpotently on
$\tau_{Q(z^{-1})}(v\otimes
t^{n})e^{-zL'_{Q(z^{-1})}(1)}(\psi^{*})^{-1}(f)$ for any $n\in {\Bbb
Z}$.)

Now we assume in addition that $f$ satisfies the $P(z)$-grading
condition. Taking $v=\omega$ in (13.49), we get
\begin{equation}
L'_{P(z)}(0)f=\psi^{*}(e^{zL'_{Q(z^{-1})}(1)}L'_{Q(z^{-1})}(0)
e^{-zL'_{Q(z^{-1})}(1)}
(\psi^{*})^{-1}(f))
\end{equation}
or equivalently, applying $e^{-zL'_{P(z)}(1)}$ and $(\psi^{*})^{-1}$
and using what we have established,
\begin{equation}
(\psi^{*})^{-1}(e^{-zL'_{P(z)}(1)}L'_{P(z)}(0)f)=L'_{Q(z^{-1})}(0)
(\psi^{*})^{-1}(e^{-zL'_{P(z)}(1)}f).
\end{equation}
If $f$ is an eigenvector of $L'_{P(z)}(0)$, then
$(\psi^{*})^{-1}(e^{-zL'_{P(z)}(1)}f)$ is an eigenvector (with the
same eigenvalue) of $L'_{Q(z^{-1})}(0)$ by (13.51). Since by
assumption $L'_{P(z)}(1)$ acts nilpotently on $f$, $e^{zL'_{P(z)}(1)}f$
is a finite sum of eigenvectors of $L'_{P(z)}(0)$ (recall (13.44)), so
that $$(\psi^{*})^{-1}(f)=
(\psi^{*})^{-1}(e^{-zL'_{P(z)}(1)}(e^{zL'_{P(z)}(1)}f))$$ is also a
finite sum of eigenvectors of $L'_{Q(z^{-1})}(0)$.  In general, since
$f$ is a finite sum of eigenvectors of $L'_{P(z)}(0)$, on each of
which $L'_{P(z)}(1)$ acts nilpotently, $(\psi^{*})^{-1}(f)$ is still a
finite sum of eigenvectors of $L'_{Q(z^{-1})}(0)$. That is,
$(\psi^{*})^{-1}(f)$ satisfies the $Q(z^{-1})$-grading condition for
$L'_{Q(z^{-1})}(0)$.

Continuing to assume that $f$ satisfies the $P(z)$-compatibility
condition and that $L'_{P(z)}(1)$ acts nilpotently on $f$, we note
that $e^{-zL'_{P(z)}(1)}f$ satisfies the same conditions. By the
discussion above, $(\psi^{*})^{-1}(e^{-zL'_{P(z)}(1)}f)$ satisfies the
$Q(z^{-1})$-compatibility condition. By (13.40) and (13.49) and the
fact that the component operators of $Y'_{Q(z^{-1})}(v, x)$ preserve
the space of elements of $(W_{1}\otimes W_{2})^{*}$ satisfying the
$Q(z^{-1})$-compatibility condition, we have, for any $v\in V$,
\begin{eqnarray}
\lefteqn{Y'_{P(z)}(v, x)f=}\nno\\
&&=\psi^{*}(e^{zL'_{Q(z^{-1})}(1)}Y'_{Q(z^{-1})}(v, x)
(\psi^{*})^{-1}(e^{-zL'_{P(z)}(1)}f))\nno\\
&&=e^{zL'_{P(z)}(1)}\psi^{*}(Y'_{Q(z^{-1})}(v, x)
(\psi^{*})^{-1}(e^{-zL'_{P(z)}(1)}f))
\end{eqnarray}
(note that the action of $e^{zL'_{P(z)}(1)}$ is defined).
Thus we have the equivalent formula
\begin{eqnarray}
(\psi^{*})^{-1}(e^{-zL'_{P(z)}(1)}
Y'_{P(z)}(v, x)f)
=Y'_{Q(z^{-1})}(v, x)
(\psi^{*})^{-1}(e^{-zL'_{P(z)}(1)}f),
\end{eqnarray}
and we know that the coefficient of each power of $x$ in $Y'_{P(z)}(v,
x)f$ satisfies the $P(z)$-compatibility condition and the
$L'_{P(z)}(1)$-nilpotency condition. Thus by induction,
we have
\begin{eqnarray}\label{1-1}
\lefteqn{(\psi^{*})^{-1}(e^{-zL'_{P(z)}(1)}
Y'_{P(z)}(v_{1}, x_{1})\cdots Y'_{P(z)}(v_{n}, x_{n})f)=}\nno\\
&&=Y'_{Q(z^{-1})}(v_{1}, x_{1})\cdots Y'_{Q(z^{-1})}(v_{n}, x_{n})
(\psi^{*})^{-1}(e^{-zL'_{P(z)}(1)}f),
\end{eqnarray}
and in particular,
\begin{eqnarray}\label{1-2}
\lefteqn{(\psi^{*})^{-1}(e^{-zL'_{P(z)}(1)}L'_{P(z)}(0)
Y'_{P(z)}(v_{1}, x_{1})\cdots Y'_{P(z)}(v_{n}, x_{n})f)=}\nno\\
&&=L'_{Q(z^{-1})}(0)Y'_{Q(z^{-1})}(v_{1}, x_{1})\cdots
Y'_{Q(z^{-1})}(v_{n}, x_{n})
(\psi^{*})^{-1}(e^{-zL'_{P(z)}(1)}f)\nno\\
&&
\end{eqnarray}
for any $n\ge 0$ and $v_{1}, \dots v_{n}\in V$.

Now we assume that $f$ satisfies both the $P(z)$-compatibility
and $P(z)$-local grading-restriction conditions.
By (13.44), $L'_{P(z)}(1)$ acts on $f$ nilpotently. {}From the
discussion above, we see that $(\psi^{*})^{-1}(f)$ satisfies the
$Q(z^{-1})$-grading condition and the formulas (\ref{1-1}) and
(\ref{1-2}) hold. Since $f$
satisfies part (b) of the $P(z)$-local grading-restriction
condition, we see {}from (\ref{1-1})--(\ref{1-2}) that the element
$$(\psi^{*})^{-1}(e^{-zL'_{P(z)}(1)}f)\in (W_{1}\otimes W_{2})^{*},$$
which satisfies the $Q(z^{-1})$-grading condition,
satisfies
part (b) of the $Q(z^{-1})$-local grading-restriction condition.
Thus
$$e^{zL'_{Q(z^{-1})}(1)}(\psi^{*})^{-1}(e^{-zL'_{P(z)}(1)}f)
=(\psi^{*})^{-1}(f)$$
also satisfies the $Q(z^{-1})$-local-grading restriction condition,
proving the ``only if'' part of the second assertion of the Theorem.

Let $f$ satisfy both the $P(z)$-compatibility condition and the
$P(z)$-local grading-restriction condition.
We already know that the components of
$Y'_{P(z)}(v, x)f$ satisfy the
$P(z)$-compatibility condition for any $v\in V$. Since $(\psi^{*})^{-1}f$
satisfies both the $Q(z^{-1})$-compatibility and
$Q(z^{-1})$-local grading-restriction conditions, the components of
$$Y'_{Q(z^{-1})}(v, x)
e^{-zL'_{Q(z^{-1})}(1)}(\psi^{*})^{-1}(f)$$
still satisfy these conditions.  Thus by (13.54) and (13.55), with
$e^{zL'_{P(z)}(1)}\circ \psi^{*}$ applied to both sides, we see that
the components of $Y'_{P(z)}(v, x)f$ satisfy  the $P(z)$-local grading
restriction-condition.

To prove the Jacobi identity (13.29), note that
$(\psi^{*})^{-1}(f)$ satisfies both the $Q(z^{-1})$-compatibility and
 $Q(z^{-1})$-local grading-restriction conditions, and thus the
Jacobi identity (13.38) for $Y'_{Q(z^{-1})}$ holds. {}From this Jacobi
identity, we obtain
\begin{eqnarray}\label{2}
\lefteqn{x_{0}^{-1}\delta
\left({\displaystyle\frac{x_{1}-x_{2}}{x_{0}}}\right)
\psi^{*}(e^{zL'_{Q(z^{-1})}(1)}
Y'_{Q(z^{-1})}(u, x_{1})
Y'_{Q(z^{-1})}(v, x_{2})(\psi^{*})^{-1}(f))}\nno\\
&&\hspace{2ex}-x_{0}^{-1} \delta
\left({\displaystyle\frac{x_{2}-x_{1}}{-x_{0}}}\right)
\psi^{*}(e^{zL'_{Q(z^{-1})}(1)}
Y'_{Q(z^{-1})}(v, x_{2})
Y'_{Q(z^{-1})}(u, x_{1})(\psi^{*})^{-1}(f))\nonumber \\
&&=x_{2}^{-1} \delta
\left({\displaystyle\frac{x_{1}-x_{0}}{x_{2}}}\right)
\psi^{*}(e^{zL'_{Q(z^{-1})}(1)}
Y'_{Q(z^{-1})}(Y(u, x_{0})v,
x_{2})(\psi^{*})^{-1}(f)).\nno\\
&&
\end{eqnarray}
Using (13.49) and (13.40) and the stability established in the last
paragraph, we can write (\ref{2}) as
\begin{eqnarray}
\lefteqn{x_{0}^{-1}\delta
\left({\displaystyle\frac{x_{1}-x_{2}}{x_{0}}}\right)
Y'_{P(z)}(u, x_{1})
Y'_{P(z)}(v, x_{2})e^{zL'_{P(z)}(1)}f}\nno\\
&&\hspace{2ex}-x_{0}^{-1} \delta
\left({\displaystyle\frac{x_{2}-x_{1}}{-x_{0}}}\right)
Y'_{P(z)}(v, x_{2})
Y'_{P(z)}(u, x_{1})e^{zL'_{P(z)}(1)}f\nonumber \\
&&=x_{2}^{-1} \delta
\left({\displaystyle\frac{x_{1}-x_{0}}{x_{2}}}\right)
Y'_{P(z)}(Y(u, x_{0})v,
x_{2})e^{zL'_{P(z)}(1)}f.
\end{eqnarray}
Thus we have proved that the Jacobi
identity holds when acting on $e^{zL'_{P(z)}(1)}f$ for any $f$ as above.
Replacing $f$ by $e^{-zL'_{P(z)}(1)}f$, which still satisfies the conditions,
we see that the Jacobi identity (13.29) holds.

It remains  only to prove the ``if'' part of the second assertion of the
Theorem. Note that in the proof above, (\ref{1-1}) and
(\ref{1-2}) were the main equalities that we had to
establish for the ``only if'' part;
the conclusion followed immediately. {}From the proof of  (\ref{1-1}) and
(\ref{1-2})
above, we see that conversely,
if $(\psi^{*})^{-1}(f)$
satisfies both the $Q(z^{-1})$-compatibility and
$L'_{Q(z^{-1})}(1)$-nilpotency conditions, we can prove
analogously, using the Jacobi identity (13.29),  that (\ref{1-1}) and
(\ref{1-2}) with $(\psi^{*})^{-1}\circ e^{-zL'_{P(z)}(1)}$ replaced
by $e^{-zL'_{Q(z^{-1})}(1)}\circ (\psi^{*})^{-1}$ hold. These two
equalities and the $Q(z^{-1})$-local grading-restriction condition for
$(\psi^{*})^{-1}(f)$ imply that
$f$ satisfies the
$P(z)$-local grading-restriction condition. \epfv

Proposition 13.3 and Theorem 13.9  give us another construction
of the $P(z)$-tensor product---the analogue of Theorem 6.3 of
\cite{HL1} (Theorem I.15 above):

\begin{theo}
The vector space  consisting of all elements of $(W_{1}\otimes W_{2})^{*}$
 satisfying the $P(z)$-compatibility
condition  and the $P(z)$-local grading-restric\-tion condition
equipped with $Y'_{P(z)}$ is
a generalized module and is equal to
the generalized module $W_{1}\hboxtr_{P(z)} W_{2}$.\epf
\end{theo}

Finally we have the analogues of the last two results in \cite{HL1}
(recalled in the introduction above).  The following result follows
immediately {from} Proposition 13.7, the theorem above and the
definition of $W_{1}\boxtimes_{P(z)} W_{2}$:

\begin{corol}
The $P(z)$-tensor product of $W_{1}$ and $W_{2}$ exists if and only if
the subspace of $(W_{1}\otimes W_{2})^{*}$ consisting of the elements
satisfying the compatibility condition and the local
grading-restriction condition, equipped with $Y'_{P(z)}$,
is a module. In this case, this module coincides with the module $W_1
\hboxtr_{P(z)} W_2$, and the contragredient module of this
module, equipped with the $P(z)$-intertwining map $\boxtimes_{P(z)}$,
 is a $P(z)$-tensor product of $W_{1}$ and
$W_{2}$,  equal to the structure $(W_{1}\boxtimes_{P(z)}
W_{2}, Y_{P(z)};
\boxtimes_{P(z)})$ constructed above.\epf
\end{corol}

{From} this result and Propositions 13.5 and 13.6, we have:

\begin{corol}
Let $V$ be a rational vertex operator algebra and $W_{1}$, $W_{2}$ two
$V$-modules. Then the $P(z)$-tensor product
$(W_{1}\boxtimes_{P(z)} W_{2}, Y_{P(z)};
\boxtimes_{P(z)})$ may be constructed as described in Corollary 13.11.\epf
\end{corol}

 {\small \sc Department of Mathematics, University of Pennsylvania,
Philadelphia, PA 19104}

{\em Current address:} Department of Mathematics, Rutgers University,
New Brunswick, NJ 08903

{\em E-mail address}: yzhuang@math.rutgers.edu

\vskip 1em

{\small \sc Department of Mathematics, Rutgers University,
New Brunswick, NJ 08903}

{\em E-mail address}: lepowsky@math.rutgers.edu

\end{document}